\documentclass[smallextended]{svjour3} 
\journalname{Mathematical Programming Computation}

\usepackage[lmargin=2cm, rmargin=2cm, tmargin=2cm, bmargin=2cm, marginpar=0pt, marginparsep=0pt]{geometry}
\makeatletter
\renewcommand{\@maketitle}{%
  \newpage
  \null
  \vspace*{-1.5cm}
  \begin{center}
    {\LARGE\bfseries \@title\par}
    \vskip 0.8em
    {\large
      \lineskip .5em%
      \begin{tabular}[t]{c}
        \@author
      \end{tabular}\par}
    \vskip 0.8em
    {\small \@date\par}
  \end{center}
  \par
  \vskip 1.0em
}
\makeatother

\usepackage{amsmath, mathtools}
\usepackage[warnings-off={mathtools-colon,mathtools-overbracket}]{unicode-math}
\usepackage{graphicx}
\usepackage[dvipsnames]{xcolor}
\usepackage[linesnumbered,ruled,vlined]{algorithm2e}
\usepackage[colorlinks,linkcolor=blue,citecolor=blue,urlcolor=magenta,linktocpage,plainpages=false]{hyperref}
\usepackage{subcaption}
\usepackage{enumitem}
\usepackage{array}
\usepackage{authblk}
\usepackage{longtable}
\usepackage{tabularx,array}
\usepackage{tablechart,booktabs,multirow,colortbl}
\usepackage{pgfplots, pgfplotstable}
\usepgfplotslibrary{fillbetween}
\usepgfplotslibrary{statistics}
\pgfplotsset{compat=1.18}
\usepackage{float}
\usepackage{tikz}
\usetikzlibrary{external}
\usetikzlibrary{arrows,shapes,plotmarks}
\usetikzlibrary{calc, matrix, patterns, patterns.meta, pgfplots.statistics, shapes.misc, positioning, arrows.meta, decorations.markings, angles}

\usepackage{academicons}
\newcommand{\orcid}[1]{\href{https://orcid.org/#1}{\textcolor[HTML]{A6CE39}{\aiOrcid}}}

\newcommand{\tagarray}{%
	\mbox{}\refstepcounter{equation}%
	$(\theequation)$%
}

\tikzset{
    cheating dash/.code args={on #1 off #2}{
        \csname tikz@addoption\endcsname{%
            \pgfgetpath\currentpath%
            \pgfprocessround{\currentpath}{\currentpath}%
            \csname pgf@decorate@parsesoftpath\endcsname{\currentpath}{\currentpath}%
            \pgfmathparse{\csname pgf@decorate@totalpathlength\endcsname-#1}\let\rest=\pgfmathresult%
            \pgfmathparse{#1+#2}\let\onoff=\pgfmathresult%
            \pgfmathparse{max(floor(\rest/\onoff), 1)}\let\nfullonoff=\pgfmathresult%
            \pgfmathparse{max((\rest-\onoff*\nfullonoff)/\nfullonoff+#2, #2)}\let\offexpand=\pgfmathresult%
            \pgfsetdash{{#1}{\offexpand}}{0pt}}%
    }
}

\newcolumntype{x}[1]{>{\centering\arraybackslash\hspace{0pt}}p{#1}}

\definecolor{msdarkblue}{RGB}{36,58,94}
\definecolor{msblue}{RGB}{0,120,215}
\definecolor{msgreen}{RGB}{16,124,16}
\definecolor{msred}{RGB}{216,59,1}
\definecolor{purple}{RGB}{128,0,128}
\definecolor{msgray}{HTML}{DFDFDF}
\definecolor{Emerald}{HTML}{00A99D}
\definecolor{RubineRed}{HTML}{ED017D}
\definecolor{pink}{HTML}{FF55A3}

\definecolor{myd}{HTML}{f0b077}
\definecolor{myp}{HTML}{ea801c}
\definecolor{mylt}{HTML}{3594cc}
\definecolor{mymt}{HTML}{8cc5e3}

\DeclarePairedDelimiter\floor{\lfloor}{\rfloor}
\newcommand{\norm}[1]{\left\lVert #1 \right\rVert}

\begin{document}

\title{Accelerating Column Generation in Highly Degenerate Integer Programming Problems with Template Pricing}
\titlerunning{Template Pricing}

\authorrunning{Marshall, Shah, and Dey}

\author{Luke Marshall \and 
        Prachi Shah \and 
        Santanu S. Dey}

\institute{
    Luke Marshall \orcid{0000-0003-0633-4004} \at
    Microsoft Research, Redmond, USA
    \email{luke.marshall@microsoft.com}     
\and
    Prachi Shah \orcid{0009-0008-7964-4875} \at
    School of Industrial and Systems Engineering, Georgia Institute of technology, Atlanta, USA
    \email{prachi.shah6795@gmail.com}     
\and
    Santanu S. Dey \orcid{0000-0003-0294-8287} \at
    School of Industrial and Systems Engineering, Georgia Institute of technology, Atlanta, USA
    \email{sdey30@gatech.edu}     
}

\date{Received: date / Accepted: date}

\maketitle
\begin{abstract}
We propose a new pricing strategy for column generation (CG), referred to as \textit{Template pricing}. This method is motivated by the desire to coordinate solutions of different pricing subproblems in order to accelerate the convergence of the CG process and simultaneously obtain good quality integer feasible solutions. Instead of finding a column with the optimal reduced cost, Template pricing tries to maximize the similarity of columns with a given template vector, while restricting the search to columns with suitable reduced cost. We present an exact and heuristic method (based on Lagrangian relaxation) to efficiently solve the Template pricing problem. We conduct extensive computational experiments on benchmark instances of the Generalized Assignment Problem (GAP). Our results demonstrate that Template pricing can significantly accelerate the CG algorithm, especially in the presence of significant degeneracy, where several benchmark GAP instances solved over 1000x faster than Dantzig pricing, and over 100x with adaptive dual-smoothing. Template pricing allows us to achieve CG optimal bounds on all 1735 ISA instances, finding stronger bounds in 43\% and improved integer solutions in 9\% of these instances than previously released.

\end{abstract}

\section{Introduction}
Column generation (CG) is used when the number of variables in a linear program (LP) is so large that it is intractable to be explicitly enumerated. Many discrete optimization problems with compact representations are reformulated using decomposition-based methods such as Dantzig-Wolfe decomposition~\cite{dantzig1960decomposition,vanderbeck2010reformulation}, to have such a large set of variables.
This is often advantageous, as the reformulated problem may provide stronger dual bounds, use a smaller working model to solve, and parallelize more easily.

While the idea behind CG is simple in theory, its effectiveness in practice is encumbered by challenges, such as efficiently generating columns of good quality, managing existing columns, maintaining the stability of the dual values, reducing the impact of degeneracy and integrating CG with effective branching and cuts~\cite{vanderbeck2010reformulation,lubbecke_column_2011,uchoa2024optimizing}. Moreover, the columns generated for different subproblems in any iteration are often uncoordinated since the pricing problems are solved independently. As a result, only a small fraction of the new columns are likely to enter the basis in the subsequent iteration, leading to slower convergence. Moreover, as Ghoniem and Sherali~\cite{ghoniem_complementary_2009} point out, the optimal restricted master problem (RMP) might not even contain columns that can be used to construct a feasible solution of the underlying integer program, let alone an optimal integer solution. We call a set of columns \textit{compatible} if they satisfy the linking constraints in the master problem and \textit{incompatible} otherwise. This is a generalization of the notion of complementarity introduced in~\cite{ghoniem_complementary_2009} for set-partitioning problems.

In this paper, we focus on coordinating the solutions of pricing subproblems to generate {nearly} compatible columns.
We propose a novel approach to generate high-quality columns using a technique that we call \textit{Template pricing}. Template pricing modifies the standard pricing problem that generates new columns. Instead of finding columns with optimal \textit{reduced cost}, Template pricing tries to maximize the similarity of columns with a given \textit{template} vector, while restricting the search to columns with suitable reduced cost. Therefore, pricing across subproblems is coordinated by choosing templates that constitute a set of compatible columns. Further, due to the altered objective, the direction of the search in the pricing problem only considers the dual values for tie-breaking, thereby reducing the impact of dual instability and primal degeneracy. Moreover, our approach can be used in conjunction with other ``dual-value stabilization'' techniques, and a combined approach may yield even greater impact.



The main contributions of this paper are the following:
\begin{itemize}
    \item We propose a new paradigm for pricing, referred to as Template pricing, a novel approach for coordinating solutions of pricing subproblems to accelerate the convergence of the CG process, while simultaneously obtaining good quality integer feasible solutions. 
    
    \item We present Template pricing in the context of set-partitioning constraints, with a primary focus on the Generalized Assignment Problem (GAP), a special subclass of such problems; however, this approach can be extended to any problem solved using CG.

    \item We present an exact approach and a Lagrangian relaxation based heuristic to solve the Template pricing problem. 

    \item We study the impact of different initialization methods and column management on the effectiveness of Template pricing, and their relationship with degeneracy.
    
    \item We conduct extensive computational experiments on benchmark instances~\cite{GAPsmall,GAPlarge,10.1007/978-3-031-26504-4_30} of GAP~\cite{oncan2007survey,cattrysse1992survey,savelsbergh1997branch,chu1997genetic,yagiura2006path}. We compare our approach with Dantzig pricing and Pessoa pricing~\cite{wentges1997weighted}. Our results demonstrate that Template pricing can significantly accelerate the CG algorithm, with some instances solving over 1000x faster. 

\end{itemize}

The rest of the paper is organized as follows. We provide background on CG in Section~\ref{sec:CGpre}. Section~\ref{sec:TemplatePricing} introduces Template pricing and compares our approach with others in the literature, which combat similar challenges. We also present a Lagrangian relaxation based approach to solve the Template pricing problem. 
In Section~\ref{sec:SolveTP} we analyze important design choices when implementing CG, such as initialization and column management of the RMP. We make informed parameter choices for all methods we evaluate, via extensive experimentation.
We present results of our computational study in Section~\ref{sec:results}, where we compare Template pricing against existing approaches, and provide evidence and insights into its superior performance. Finally, we conclude and discuss future research in Section~\ref{sec:Conclusion}.

\section{Column Generation (CG) Preliminaries}\label{sec:CGpre}

In this section, we introduce CG for problems where the linking constraints are of the set-partitioning type. This constitutes a large class of problems solved using CG and covers GAP~\cite{savelsbergh1997branch,pigatti2005stabilized}, vehicle routing~\cite{desaulniers2006column,fukasawa2016branch,baldacci2008exact,pecin2017improved}, graph coloring~\cite{mehrotra1996column,held2012maximum}, bin packing~\cite{valerio1999exact}, and a wide range of scheduling problems such as crew scheduling~\cite{haase2001simultaneous,desaulniers1997crew,barnhart1998branch} and the gate assignment problem~\cite{li2021using}. The compact formulation is represented as, 

\begin{center}
\begin{subequations}\label{eq:compact_mip}
\renewcommand*{\arraystretch}{1.5}
\begin{tabularx}{\textwidth}{>{\raggedleft\arraybackslash}Xrclll>{\raggedleft\arraybackslash}X}
    $\displaystyle \min_{\mathbf{x}}$ & \multicolumn{4}{l}{
        $\quad\displaystyle \sum_{i\in I} \mathbf{c}^\top_{i}\mathbf{x}_{i} $} & & \tagarray\label{eq:vectorcost}\\


    s.t. & $\quad\displaystyle \sum_{i\in I} \mathbf{e}_j^\top \mathbf{x}_{i}$ & $=$ & $\mathbf{1}$ & $\qquad \forall j\in J$, & & \tagarray\label{eq:partition} \\

    & $\quad\displaystyle \mathbf{A}_{i} \mathbf{x}_{i}$ & $\leq$ & $\mathbf{b}_{i}$ & $\qquad\forall i\in I$, & & \tagarray\label{eq:block-constraints}\\

    & $\quad \mathbf{x}_{i}$ & $\in$ & $\{0,1\}^{|J|}$ & $\qquad\forall i\in I$, &&\tagarray\label{eq:integral}
\end{tabularx}
\end{subequations}
\end{center}
where $e_j$ is the unit vector in the $j^{\text{\,th}}$ direction. Since GAP is a prototypical set-partitioning problem, we discuss \eqref{eq:compact_mip} in this context. Here, the variable $\mathbf{x}_i \coloneqq (x_{i1},\; x_{i2},\; \dots,\; x_{i|J|})^\top$ is a vector with binary elements $x_{ij}$ that represents if job $j\in J$ is assigned to machine $i\in I$. The set-partitioning constraint (\ref{eq:partition}) ensures that all jobs are assigned exactly once. The complicating constraints (\ref{eq:block-constraints}) are a single knapsack constraint for each machine that enforces that the assignment of job demands do not exceed machine capacity.  
Henceforth in the paper, it will often be convenient to refer to the $i$ index as machine and $j$ index as job even for the general set-partitioning instances. 

This formulation can be decomposed using Dantzig-Wolfe reformulation, and by taking its linear relaxation gives the following formulation: 
\begin{center}
    \begin{subequations}\label{eq:RMP}
    \renewcommand*{\arraystretch}{1.5}
    \begin{tabularx}{\textwidth}{>{\raggedleft\arraybackslash}Xrclll>{\raggedleft\arraybackslash}X}
		$\displaystyle \min_{\lambda}$ & \multicolumn{4}{l}{
			$\quad\displaystyle \sum_{i \in I} \sum_{\mathbf{v}_{ip} \in P_i} \left(\mathbf{c}_i^\top\mathbf{v}_{ip}\right) \lambda_{ip} $} && \tagarray\\
		
		s.t. & $\quad\displaystyle \sum_{i \in I} \sum_{\mathbf{v}_{ip} \in P_i} \mathbb{1}_{\{\mathbf{e}_j^\top\mathbf{v}_{ip} = 1\}} \lambda_{ip}$ & $=$ & $1$ & $\quad\forall j\in J $ & $\quad (\pi_j)$ & \tagarray\label{eq:mpcon1} \\
		
		& $\quad\displaystyle \sum_{\mathbf{v}_{ip} \in P_i} \lambda_{ip}$ & $=$ & $1$ & $\quad\forall i\in I$ & $\quad (\mu_i)$ & \tagarray\label{eq:mpcon2}\\
		
		& $\quad\lambda_{ip}$ & $\geq$ & $ 0$ & $\quad\forall i\in I,\; \mathbf{v}_{ip}\in P_i.$ && \tagarray 
    \end{tabularx}
    \end{subequations}\label{eq:cg}
\end{center}
where $P_i = \left\{\mathbf{x} \in \left\{0, 1\right\}^{|J|} \mid \mathbf{A}_i \mathbf{x} \leq \mathbf{b}_i\right\}$, and $\lambda_{ip}$ is the variable corresponding to the column $\mathbf{v}_{ip} \in P_i.$ The notation $\mathbb{1}_{\{\mathbf{e}_j^\top\mathbf{v}_{ip} = 1\}}$ represents the value $1$ if the $j^{\text{\,th}}$ job is assigned to the $i^{\text{\,th}}$ machine in column $\mathbf{v}_{ip}$ and $0$ otherwise. 
The dual variables for \eqref{eq:mpcon1} and \eqref{eq:mpcon2} are $\pi_i$ and $\mu_j$, respectively. Enumerating all vectors in $P_i$ is often intractable, so to solve \eqref{eq:cg}, a restriction on the subset of vectors $V_i \subseteq P_i$ is instead solved iteratively. The model defined on $V_i$ is called the Restricted Master Problem (RMP). 
At every iteration, the pricing problems are solved for each subproblem $i$ to identify the columns to add to $V_i$, or to prove optimality (and terminate) if none exist. For the column $\mathbf{x}\in P_i$, to be added to $V_i$, it must have a negative reduced cost, also referred to as a good reduced cost. Therefore, it must satisfy the following, for some small tolerance $\varepsilon>0$,
\begin{equation}
\sum_{j \in J} \left(c_{ij} - \pi_j\right) x_j \leq \mu_i - \varepsilon. \label{eq:goodredcost}
\end{equation}
At each iteration, at least one column having good reduced cost must be added to the RMP across all machines to guarantee convergence. If no such columns exist for all $i$, this certifies that the solution of RMP is optimal for the original master problem \eqref{eq:cg}.
Formally, at each iteration, having corresponding dual values $(\pi, \mu)$, the pricing subproblem chooses columns from the set $S^i_{(\pi,\mu)}=\left\{ \mathbf{x}\in P_i \;\middle|\; \eqref{eq:goodredcost} \right\}$ or shows that $S^i_{(\pi,\mu)}=\emptyset$ for all $i\in I$.

\subsection{Relation with Lagrangian Relaxation}\label{sec:LR}

Though not the focus of this paper, Lagrangian Relaxation (LR) is closely related to CG. We include a discussion of LR to provide context for the pricing strategies explored in Section \ref{sec:smoothing} (dual smoothing), and also because prior studies report strong practical performance on GAP under degeneracy~\cite{savelsbergh1997branch}. 

Both LR and CG are decomposition techniques driven by dual information~\cite{geoffrion2009lagrangean,fisher1981lagrangian}. LR relaxes hard constraints by incorporating them into the objective function with associated penalties (Lagrange multipliers). We outline the LR derivation in the remainder of this section. For GAP, relaxing the set partitioning constraints \eqref{eq:partition} in LR yields the same bounds as the CG RMP. By relaxing these linking constraints, LR decomposes GAP into independent binary knapsack problems -- again, similar to the CG decomposition. Relaxing the set partitioning constraints transforms the problem into:
\begin{align*}
\max_{\pi} \left\{ \min_{\mathbf{x}} \left\{ \sum_{i \in I} \mathbf{c}^\top_{i}\mathbf{x}_{i} + \sum_{j\in J} \pi_j \left(1 - \sum_{i\in I} \mathbf{e}_j^\top \mathbf{x}_{i} \right) \;\middle|\; \eqref{eq:block-constraints}, \eqref{eq:integral} \right\} \right\} &\\
=\max_{\pi} \left\{ \sum_{j\in J} \pi_j + \sum_{i \in I} \min_{\mathbf{x}} \left\{ \sum_{j \in J} \left( c_{ij} - \pi_j \right) x_{ij} \;\middle|\; \eqref{eq:block-constraints}, \eqref{eq:integral} \right\} \right\} &.
\end{align*}

 Moreover, these Lagrange multipliers correspond directly to the dual values of the CG RMP.  LR explicitly updates the multipliers iteratively (via subgradient, bundle, level methods etc.) depending on the chosen implementation. LR uses the gradient $\nabla\pi$ to update the multipliers at each step: 
$$\nabla\pi_j = 1 - \sum_{i \in I} x_{ij} = 1 - \sum_{j\in J}
    \mathbb{1}_{\{\mathbf{e}_j^\top\mathbf{v}_{ip} = 1\}}.$$

In degenerate GAP instances, LR can often outperform CG since it avoids solving an LP at each iteration. Nevertheless, LR has its own drawbacks: no primal feasibility guarantees (often needs additional heuristics to obtain primal solutions), sensitive to step sizes, potential oscillations and slow convergence, limited coordination across subproblems, and often terminates at a near-optimal solution (without clean certification). 

Note that, similar to \eqref{eq:cg}, subproblems in LR yield solutions of the form of $v_{ip}$. 
Steps are performed iteratively until some termination criteria; for example, when the change in objective value between two iterations is less than a given threshold.


\subsection{Column Generation Pricing Strategies}

The most natural approach for choosing columns from $S^i_{(\pi,\mu)}$ is \textit{Dantzig} pricing. It uses \eqref{eq:goodredcost} as the objective function and checks whether the resulting minimal value has good reduced cost, specifically: 
\begin{equation}\label{eq:dantzig-pricing}
   rc(i):= \min_{x\in P_i} \left\{ \sum_{j\in J} \left(c_{ij} - \pi_j\right) x_j \right\} - \mu_i \leq  - \varepsilon.
\end{equation}
It is not necessary to find columns with optimal reduced cost, in fact, it is often unlikely that columns with minimum reduced cost are the best choice to reduce time and iterations for the RMP to converge. Other objective functions, also known as \textit{pricing strategies}, for choosing solutions from $S^i_{(\pi,\mu)}$ have been proposed in the literature in the context of column generation~\cite{lubbecke_column_2011} and more generally for simplex method, such as deepest-cut in the dual space~\cite{vanderbeck1994}, and steepest-edge~\cite{forrest1992steepest,goldfarb1977practicable} and its variants as Devex~\cite{Harris1973PivotSM} and Lambda~\cite{Bixby1992VeryLL} pricing. Although CG can be considered an extension of primal simplex, many pricing strategies applicable for primal simplex are not practical for CG due to non-linearity, since CG pricing involves solving an optimization problem, it is computationally preferable to have a linear objective function.

\textcolor{black}{More generally, an important unresolved question is whether there exists a pivoting rule under which the simplex algorithm operates in polynomial time. See for example, the smoothened analysis~\cite{spielman2009smoothed,dadush2018friendly} and the paper~\cite{black2025exponential} and references therein for a comprehensive discussion of proposed rules and several lines of research.} 

\subsubsection{Disadvantages with Dantzig Pricing}

Dantzig pricing can be an excellent choice for generating columns. It simplifies $S^i_{(\pi,\mu)}$ by moving \eqref{eq:goodredcost} into the objective; it yields valid bounds on the RMP at each iteration; and is typically the way to guarantee termination at optimality. However, Dantzig pricing suffers on instances that have degeneracy. 

Degeneracy is known to be associated with unstable dual values in the RMP when solving with the simplex algorithm. This is unfortunate, since primal simplex is a good choice for solving the RMP in CG, due to the computational efficiency of warm starts when adding columns. 
Unstable dual values provide bad search directions for pricing and
fluctuate greatly between CG iterations. This issue is well known in the literature and there exist many papers on dual stabilization techniques to improve convergence in CG; see for example~\cite{wentges1997weighted,du1999stabilized,elhallaoui2005dynamic,ben2006dual,amor2004stabilization,rousseau2007interior,briant2008comparison,amor2009choice,ghoniem_complementary_2009,elhallaoui2010multi,benchimol2012stabilized,pessoa2013out,gschwind2016dual,pessoa2018automation,costa2022stabilized}.
These approaches encompass a variety of strategies, including trust-region methods, interior-point stabilization, constraint aggregation and coordination schemes, as well as dual-optimal inequalities.

Dual value instability also yields issues with uncoordinated subproblems. There are computational advantages to solving subproblems independently, such as parallelization. However, the columns generated for different subproblems may be incompatible, that is, may violate the RMP constraints and therefore impossible to simultaneously add to the RMP. Dantzig pricing implicitly uses dual values for coordination, but it is fragile, especially with unstable dual values, since it relies on potentially small differences in objective coefficients.  In particular, observe that for job $j$, the objective coefficient in \eqref{eq:dantzig-pricing} relies on the dual value $\pi_j$, which is shared across every subproblem for each machine and perturbed by the $c_{ij}$ costs. Jobs corresponding to large dual values that outweigh $c_{ij}$ perturbations are likely to be selected by several machines, whereas jobs corresponding to small dual values might not be selected at all. The generated columns often compete for the same jobs and are unlikely to enter the RMP basis simultaneously in the subsequent iteration. This wasted computational effort leads to less improvement in the RMP solution across iterations and therefore has slower convergence. 

To avoid the issue of incompatible columns, Ghoniem and Sherali \cite{ghoniem_complementary_2009} sequentially solve the pricing problem, while cumulatively adding constraints that enforce compatibility within the sequence. This ensures the compatibility of the columns and yields an integer feasible solution; however, the quality of solution is dependent on the choice of sequence ordering and is more computationally expensive, since it prohibits parallel computation. 
Alternatively, there are dynamic constraint aggregation (DCA) methods,  see for example~\cite{elhallaoui2005dynamic,elhallaoui2010multi,benchimol2012stabilized,bouarab_dynamic_2017}, that group jobs together and modifies the pricing such that it first tries to find columns satisfying this grouping while having good reduced cost. In this way, it implicitly yields nearly compatible columns at each iteration.  If no suitable columns are found, the aggregation is refined.  This approach provides dual stabilization, parallelization, and uses a smaller RMP; at the cost of increased complexity.  Specifically, the refinement step may require solving another CG problem, which may have its own issues with degeneracy.

\subsubsection{Pricing with Dual Smoothing}
\label{sec:smoothing}

Among the numerous approaches to dealing with dual stabilization discussed above, our paper focuses specifically on the impact of pricing strategies. We exclude approaches that either modify the RMP or use alternative algorithms to solve the RMP, so that we can isolate the contribution of pricing.


In this section, we describe a strong dual smoothing technique that is used as a baseline in our study.
At a high level, dual smoothing constructs a smoothed dual vector $\tilde{\pi}_t$ for each iteration $t$, using information gathered across iterations.  It then uses $\tilde{\pi}_t$ rather than the raw RMP duals $\pi_t$ to define the pricing objectives. Two classical forms are:
\begin{itemize}
    \item \textit{Neame}~\cite{neame_nonsmooth_2000}: a static convex combination of the previous smoothed duals and the current duals, $$\tilde{\pi}_t = \alpha_0\tilde{\pi}_{t-1}+\left(1-\alpha_0\right)\pi_t$$

    \item \textit{Wentges}~\cite{wentges1997weighted}: a static convex combination of the most recent improving dual vector $\hat{\pi}_{t-1}$ (from the most recent iteration that improved the RMP objective) and the current duals, $$\tilde{\pi}_t = \alpha_0\hat{\pi}_{t-1}+\left(1-\alpha_0\right)\pi_t$$
\end{itemize}
with the static mixing parameter $\alpha_0 \in [0,1]$. In our baseline implementation we use the advanced smoothing techniques from \textit{Pessoa}~\cite{pessoa_automation_2018}. 

Pessoa generalizes Wentges by adding (i) an adaptive mixing parameter and (ii) a directional component derived from subgradient information; coupled with a limited backtracking procedure designed to ensure that pricing still produces at least one column with negative reduced cost (when possible).  This backtracking progressively reduces the impact of smoothing and ultimately reverts Pessoa to standard Dantzig pricing.

In our GAP setting, Pessoa's directional component uses the subgradient information from Lagrangian Relaxation, defined in Section~\ref{sec:LR}.  The elements in the subgradient correspond to each job, and reflects whether the solutions found by independently solving each subproblem (per machine) over-cover or under-cover that job in the current iteration. Pessoa uses this information to penalize over-covered jobs and prioritize under-covered jobs, providing directional correction beyond simple convex smoothing. 

While highly effective, the mathematics behind Pessoa is substantially more complex than Wentges, Neame, and even Template pricing.  For completeness, we provide the full details of our implemented Pessoa variant in Appendix~\ref{sec:pessoa}.

\section{Template Pricing} \label{sec:TemplatePricing}

Template pricing aims to find a set of new columns (one for each subproblem) that maximizes similarity to a corresponding set of template vectors, while also having good reduced cost. These template vectors coordinate the subproblems, even though each subproblem is solved independently. Templates can be chosen from existing columns or generated by some heuristic and serve as a reference point to guide the search.

Formally, suppose we are given $\left\{y^i\right\}_{i\in I}$, a near-optimal solution to the relaxed RMP that satisfies the linking constraints. We wish to guide the search to converge towards this solution. We refer to these $y^i\in \mathbb{R}^{|J|}$ vectors as \textit{templates} and define Template pricing as the strategy that finds columns from $S^i_{(\pi,\mu)}$ with maximum similarity to the corresponding template; specifically: 
\begin{equation}
\min_{x\in S^i_{(\pi,\mu)}} \, \left(-d\left(x,y^i\right),\quad \sum_{j\in J}\left(c_{ij} - \pi_j\right)x_j\right),\label{eq:template-pricing}
\end{equation}
a bi-objective optimization, which can be interpreted as:
\begin{eqnarray}
&\min& \sum_{j\in J}\left(c_{ij} - \pi_j\right)x_j \label{eq:objouter}\\
&\textup{s.t.}& x \in {\arg \max}_{x\in S^i_{(\pi,\mu)}} \left\{ d(x,y^i)\right\}. \label{eq:objinner}
\end{eqnarray}
Template pricing \eqref{eq:template-pricing} defines a hierarchical optimization problem where the primary objective is to maximize similarity to the template, with ties broken using the reduced costs. The domain $x\in S^i_{(\pi,\mu)}$ ensures that the optimization produces columns with negative reduced cost. Equivalently, we can view this as the inner problem~\eqref{eq:objinner} selecting a column that maximizes similarity to the template, and then the outer problem~\eqref{eq:objouter} acting as a tie-breaker, choosing among the optimal solutions of \eqref{eq:objinner} the columns with minimum reduced cost. 

\subsection{Designing Templates}
We propose designing templates in a \textit{static} or \textit{dynamic} fashion, depending on whether they change during the CG process. For example, a static template could be constructed from a heuristic solution to \eqref{eq:compact_mip}, or the solution $x^*$ of its LP relaxation, that is, $y^i = \{ x^*_{ij} \}_{j\in J}$ for $i\in I$.
Other ideas include using convex combinations of near-optimal linear relaxation solutions, or Lagrangian lower bound solutions.


We propose dynamic templates based on the incumbent RMP solution, updated periodically. Let $\lambda^*$ be the solution of the RMP at some iteration. The template for each machine $i$ is defined as: 
\begin{eqnarray}\label{eq:dynamic}
y^i = \sum_{v_{ip} \in P_i} \lambda_{v_{ip}}^* v_{ip}, 
\end{eqnarray}
that is, the template is the projected primal solution of the RMP.

Initial experiments showed that static templates are beneficial, particularly when chosen well; however, the dynamic template \eqref{eq:dynamic} consistently produced improvements. Therefore, we present only results using \eqref{eq:dynamic} in the paper. 



\subsection{Designing Similarity Function}
The function $d$ measures the similarity of $x$ to the template $y^i$, where larger values indicate greater similarity. 
In this paper, we define $d$ as follows: 

\begin{equation} \label{eq:distance_fn}
    d(x, y) = \sum_{j\in J} f\left(y_j\right) x_j, \text{ where } f(y_j) = \left\{ \begin{array}{cl}
     \phantom{-}1& \quad\textup{if } y_j \in \left(1- \delta,\, 1\right] \\
    \phantom{-}0& \quad\textup{if } y_j \in \left[\delta,\, 1-\delta\right]\\
    -1& \quad\textup{if } y_j \in \left[0,\, \delta\right).
    \end{array}\right. ,
\end{equation}
for some $\delta \in [0, 0.5]$; we choose $\delta=10^{-6}$.
The intuition behind (\ref{eq:distance_fn}) is as follows: we wish to set $x_j=1$ when $y_j\approx 1$ and $x_j=0$ when $y_j\approx 0$.  When $y_j$ is fractional (neither close to 0 or 1) we wish to choose $x_j$ by tie breaking on the reduced cost.
Specifically, we choose $x_j \in \{0,1\}$ such that $f(y_j)x_j$ is maximized, so $f(y_j\approx 1) = 1$ encourages $x_j=1$ and $f(y_j\approx 0) = -1$ encourages $x_j=0$. When $y_j$ is fractional, $f(y_j)=0$ does not encourage the choice of $x_j$ in either direction. In other words, maximizing $d(x,y)$ favors ``soft fixing" variables to their template values, whenever they are near integral.


Templates that are highly fractional with a small choice of $\delta$, lead to $d(x,y^i) \approx 0$, that is, Template pricing reverts to Dantzig pricing.  It could be quite reasonable to combine Template pricing with dual smoothing, like Pessoa, by adjusting the tie-breaking objective; for simplicity, however, we do not explore this in our paper.

\subsection{Advantages of Template Pricing}
We list below several advantages of Template pricing.
\begin{enumerate}[itemsep=0.5em]
    \item \textit{Coordinated independent subproblems:} As discussed above, solving pricing subproblems independently offers computational advantages such as parallelism, but it is likely to produce incompatible columns. Template pricing coordinates these independent subproblems via compatible templates. By attempting to find columns similar to the templates, they are more likely to be compatible. Thus, Template pricing supports parallelization while reducing column incompatibility.
   
    \item \textit{Column stabilization:} 
    Other pricing algorithms attempt to stabilize dual values to reduce fluctuations between iterations and improve convergence.  Template pricing attempts to stabilize the generated columns themselves, such that the columns from one iteration to the next are similar. Section~\ref{sec:results} demonstrates that Template pricing reduces the number of simplex pivots per iteration in the RMP. Each pivot is also faster, which significantly accelerates the overall convergence speed.

    \item \textit{Integrality:}
    Using \eqref{eq:dynamic} for our template and \eqref{eq:distance_fn} for similarity in our Template pricing ensures compatibility and reinforces column stabilization. This strategy also encourages the RMP towards integer feasible solutions throughout the CG process. In addition to feasible solutions, the resulting integral solution yields a smaller basis, which further reduces pivot times.
    
    \item \textit{Degeneracy:} Prevalent in partitioning based CG problems, particularly when the job/machine ratio is high~\cite{savelsbergh1997branch}, degeneracy leads to multiple dual optimal solutions, making the direction indicated by a single dual solution less reliable. This leads to degenerate pivots and wasted computational effort.
    Template pricing reduces the impact of degeneracy by instead focusing on template similarity; specifically, these dual values have less control in the generation of columns (provided the reduced costs are negative). 
    In this regard, Template pricing shares similarities with DCA methods~\cite{elhallaoui2005dynamic,elhallaoui2010multi,bouarab_dynamic_2017} that iteratively fix a subset of jobs to some machines to create a reduced non-degenerate problem. Likewise, Template pricing prefers anchoring the jobs in the template on a machine, but does not enforce this strictly. 
    
   \item \textit{Simplicity:}  
   Template pricing is easy to implement (see Appendix~\ref{sec:hie}), since state-of-the-art solvers like Gurobi and HiGHS have good implementations for solving hierarchical optimization problems exactly.
   As an alternative, the next section presents a heuristic implementation of Template pricing based on Lagrangian relaxation.
In comparison,  other rules for pricing such as lambda pricing~\cite{bixby1992very} and DCA~\cite{elhallaoui2005dynamic,elhallaoui2010multi,bouarab_dynamic_2017} rules are significantly more challenging to implement.

\end{enumerate}


\subsection{Lagrangian-relaxation based approximation for Template pricing}


As discussed in the previous section, most state-of-the-art solvers provide built-in support for hierarchical objectives. However, one drawback of hierarchical optimization is that during the tie-breaking phase one must introduce an additional constraint to the feasible region of the pricing problem (see Appendix~\ref{sec:hie}). This modification can reduce the efficiency of specialized solution methods for the pricing problem, such as dynamic programming (DP) algorithms for knapsack problems. To avoid this issue, instead of solving the hierarchical optimization problem exactly, we investigate an approximate approach based on Lagrangian relaxation combined with bisection search. The main advantage of this method is that it preserves the original feasible region of the pricing problem. In particular, our approach uses a DP algorithm for solving the binary knapsack problem, adapted from the implementation in  SCIP~\cite{SCIPOptSuite10}.

We begin by verifying that $S^i_{(\pi,\mu)} \neq \emptyset$, that is, solve the Dantzig pricing problem and check for reduced cost.
If \eqref{eq:dantzig-pricing} does not hold, we terminate since $S^i_{(\pi,\mu)} = \emptyset$. This step is not strictly necessary, but it is extremely fast in our setting, and the optimal reduced cost can be used to improve bounds and terminate early.
Assume that $S^i_{(\pi,\mu)} \neq \emptyset$, we proceed as follows. At each iteration of our algorithm we solve the following Lagrangian relaxation problem corresponding to the ``Lagrangian parameter" $\alpha\geq 0$:
\begin{equation}\label{eq:lagk}
\mathrm{OPT}(\alpha)
:=
\min_{x\in P_i} \left\{
 -d(x,y^i)
 + \alpha \sum_{j\in J} (c_{ij}-\pi_j)\,x_j
\right\}.
\end{equation}
Let $x^\alpha$ denote an optimal solution of $\mathrm{OPT}(\alpha)$. Our goal is to determine the smallest value of $\alpha$ for which $x^\alpha \in S^i_{(\pi,\mu)}$. Intuitively, this corresponds to maximizing the similarity of the generated column to the template, while ensuring that it maintains a negative reduced cost. 

To initialize, we set $l := 0$, $u := +\infty$,
and $\alpha := \alpha_{\mathrm{init}} > 0$. At each iteration we solve (\ref{eq:lagk}) and then update the interval $[l,u]$ as follows. If $x^\alpha \in S^i_{(\pi,\mu)}$, set $u := \alpha$; otherwise set $l := \alpha$. In the former case, $\alpha$ is already sufficiently large to yield a solution in $S^i_{(\pi,\mu)}$, whereas in the latter case it is too small.

The value of $\alpha$ for the next iteration is selected via bisection search as
\begin{eqnarray*}
\alpha := \begin{cases}
\dfrac{l+u}{2}, & \text{if } u < +\infty,\\
2\alpha, & \text{if } u = +\infty.
\end{cases}
\end{eqnarray*}
The procedure terminates when $\frac{u-l}{l} \leq 10^{-3}$, or if it can prove optimality for a solution with negative reduced cost, that is, using $\floor{\alpha\mu_i-\text{OPT}(\alpha)}$ as a valid upper bound for maximum template similarity. We return $x^\alpha$ corresponding to $\alpha = u$.
In our implementation, we set $\alpha_{\mathrm{init}} = 0.5$ when solving the pricing problem for a given machine for the first time. For subsequent pricing problems associated with the same machine, $\alpha_{\mathrm{init}}$ is set equal to the optimal value of $\alpha$ obtained in the preceding CG iteration.

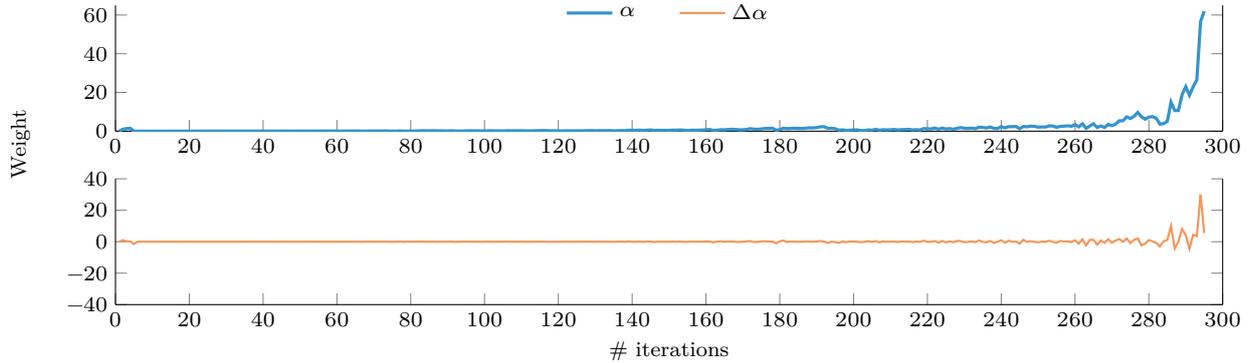
\begin{figure}[htbp]
\centering
\begin{tikzpicture}
\begin{axis}[
  name=left fig,
  axis x line*=bottom,
  xmin=0, xmax=300,
  ymin=0, ymax=65,
  ylabel={Weight},
  ylabel near ticks,
  ylabel style={at={(axis description cs:-0.07,-0.1)}},
  tick label style={font=\small}, 
  width=0.95\columnwidth,
  height=3.25cm,
  legend cell align={left},
  legend style={/tikz/every even column/.append style={column sep=0.5cm}, font = \small, at={(0.5,1.1)},legend cell align={left}, anchor=north,legend columns=-1, draw=none, fill=none}
] 
    \addplot [draw=mylt, very thick] table[x=iteration, y=a2] {c10400_LT_weights_all.tsv};
    \addlegendentry{$\alpha$}
    \addlegendimage{Peach, thick} \addlegendentry{$\Delta \alpha$}
\end{axis}

\begin{axis}[
  name=left fig,
  axis x line*=bottom,
  xmin=0, xmax=300,
  ymin=-40, ymax=40,
  xlabel={\# iterations},
  tick label style={font=\small}, 
  width=0.95\columnwidth,
  height=3.25cm,
  legend cell align={left},
  shift={(0, -2.3cm)},
  legend style={/tikz/every even column/.append style={column sep=0.5cm}, font = \small, at={(0.5,1.1)},legend cell align={left}, anchor=north,legend columns=-1, draw=none, fill=none}  
] 
    \addplot [Peach, thick] table[x=iteration, y=delta] {c10400_LT_weights_all.tsv};
\end{axis}
\end{tikzpicture} 
\caption{The optimal $\alpha$ weight and delta change per iteration for a single machine in LT pricing for the instance c10400.}
\label{fig:c10400-delta-weight}
\end{figure}

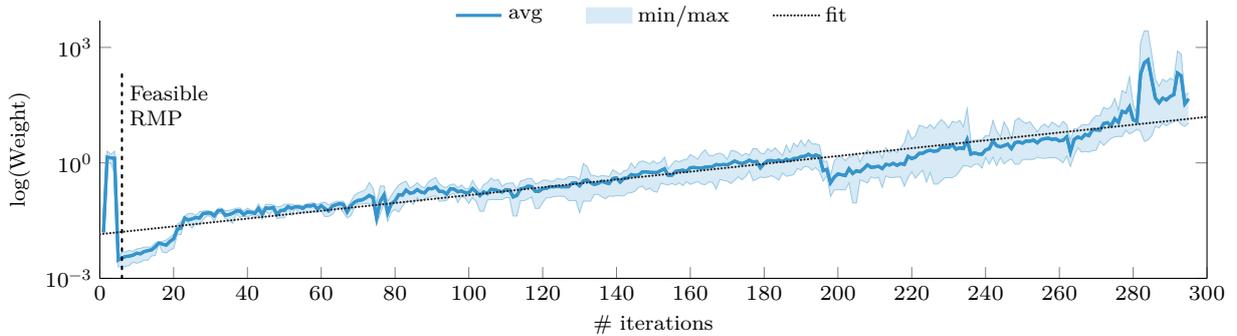
\begin{figure}[htbp]
\centering
\begin{tikzpicture}
\begin{axis}[
  name=left fig,
  axis x line*=bottom,
  xmin=0, xmax=300,
  ymin=0.001, ymax=5000,
  ymode=log,
  ylabel={log(Weight)},
  ylabel near ticks,
  xlabel={\# iterations},
  tick label style={font=\small}, 
  width=0.95\columnwidth,
  height=5cm,
  legend cell align={left},
  legend style={/tikz/every even column/.append style={column sep=0.5cm}, font = \small, at={(0.5,1.1)},legend cell align={left}, anchor=north,legend columns=5, draw=none, fill=none}
] 

\addplot [draw=mylt, very thick] table[x=iteration, y=average] {c10400_LT_weights_all.tsv};

\addplot [name path=lb, draw=mylt!50, forget plot, ultra thin] table[x=iteration,y=min] {c10400_LT_weights_all.tsv};
\addplot [name path=ub, draw=mylt!50, forget plot, ultra thin] table[x=iteration,y=max] {c10400_LT_weights_all.tsv};
\addplot[mylt!20, draw=none] fill between[of=lb and ub];

\addplot [draw=mylt, very thick, forget plot] table[x=iteration, y=average] {c10400_LT_weights_all.tsv};

\addplot [draw=Black, thick, cheating dash=on 0.0pt off 1.5\pgflinewidth, line cap=round, line join=round, domain=0:310, samples=201] (x,{0.014*exp(0.0234*x)});

\draw [draw=black, cheating dash=on 1.0pt off 3\pgflinewidth, line cap=round, line join=round, thick] (axis cs: 6, 0.001) -- (axis cs: 6, 200);

\node[anchor=west, text width=1.5cm] at (axis cs:6,30) {\footnotesize \textcolor{black}{{Feasible RMP}}};

\legend{avg, min/max, fit}
\end{axis}
\end{tikzpicture} 
\caption{Average of the optimal $\alpha$ weights across all machines per iteration in LT pricing for instance c10400. Min/Max values of $\alpha$ across machine is shown by shaded region. The fit plotted above corresponds to: $0.014 e^{0.0234x}$}
\label{fig:c10400-weights-all}
\end{figure}

We call the above procedure the Lagrange Template (LT) pricing method.
Figure~\ref{fig:c10400-delta-weight} shows the optimal $\alpha$ obtained via LT pricing at each CG iteration, together with $\Delta\alpha$ (the change in the
optimal $\alpha$ between successive CG iterations), for a single machine on one instance from our test library. Observe that optimal $\alpha$ increases significantly as CG process approaches termination. Larger values of $\alpha$ place greater weight on the reduced-cost term when solving the pricing problem, indicating that near termination it becomes increasingly difficult to identify columns with sufficient negative reduced cost; consequently, similarity to the template becomes less critical. The $\Delta\alpha$ remains consistently small (except near CG termination), suggesting that warm-starting LT with the previous value of $\alpha$ is beneficial.
Figure~\ref{fig:c10400-weights-all} reports the average optimal $\alpha$ across all machines, along with the minimum and maximum values of $\alpha$ over machines at each CG iteration. Note that the vertical axis is on a logarithmic scale, and the growth of optimal $\alpha$ appears to follow an exponential trend.

\section{Key Algorithmic Considerations}\label{sec:SolveTP}
To make the experimental analysis as fair as possible, this section covers some important considerations for parameter selection and experimental design.  In particular, we discuss the benchmark instances and show some initial results that guide our choices.

\subsection{Instances}

GAP has been extensively studied with classic sets of instances dating back at least 1994~\cite{cattrysse_set_1994}. Those 60 instances are now trivial and can be solved within 1 second; however, some instances generated in 1997~\cite{chu_genetic_1997} and extended in 2006~\cite{yagiura2006path} are still computationally challenging. At the time of writing, we have not seen evidence that all instances have been solved to optimality; however, they have all been solved within less than 0.1\% optimality gap.  We focus the core of our experiments on these 57 instances~\cite{chu_genetic_1997,yagiura2006path}. More recently, however, in 2023 \cite{10.1007/978-3-031-26504-4_30}, the classic set of instances have been extended even further by performing an instance space analysis, which attempts to systematically design a set of instances that cover a range of various metrics. We also evaluate these additional 1735 instances and provide improved bounds.  The aim for \cite{10.1007/978-3-031-26504-4_30} was to generate a diverse set of instances and evaluate with fast heuristics.  Their computational study evaluated the solutions with a time limit of 15 minutes, whereas we focus on the convergence of CG at the root node. We refer to the 57 instances as \textit{Yagiura}~\cite{yagiura2006path} and the 1735 instances as \textit{ISA}~\cite{10.1007/978-3-031-26504-4_30}.

\subsection{Degeneracy}

Savelsbergh~\cite{savelsbergh1997branch} identified that CG with Dantzig pricing has slow convergence when the ratio of jobs to machines is greater than five, and noted that LR is typically more efficient for these harder instances. The job to machine ratio is an excellent proxy for measuring degeneracy in GAP, that is, large ratios are correlated to high degeneracy and unstable dual values, which is particularly challenging for Dantzig pricing. Lübbecke and Desrosiers \cite{lubbecke_selected_2005} suggest that instances having a ratio over ten are considered to have massive degeneracy.  The instances we consider in our experiments have a minimum ratio of 5, average of $28 \pm 18$, and maximum of 80. They are excellent examples to highlight the impact of degeneracy. Our results in Section~\ref{sec:results} show that CG with Template pricing can be competitive against LR, even with high degeneracy, indicating there are additional factors to consider when choosing between CG and LR.

\subsection{Algorithms}\label{sec:algorithms}

Our study aims to evaluate the impact of various pricing strategies. We refrain from exploring algorithms that explicitly and iteratively modify the RMP or solving the RMP with interior point methods to implicitly provide dual stabilization. However, we consider one minor change to the RMP that is commonly applied in CG for dual stabilization; we replace set partitioning with cover constraints. With this in mind, our evaluation considers the following pricing strategies: 
\begin{description}
    \item[$(D)$] Dantzig -- no stabilization,
    \item[$(P)$] Pessoa -- directional and adaptive smoothing,
    \item[$(LT)$] Lagrange Template -- heuristic template, and
    \item[$(MT)$] MIP Template -- exact template.
\end{description}

As a final exception, we also include a limited comparison with Lagrangian Relaxation $(LR)$ -- due to Savelsbergh's guidance~\cite{savelsbergh1997branch} on degenerate instances. Specifically, on Yagiura we compare $\{D,P,LT,MT\}$, and ISA $\{P,LT,LR\}$. We exclude D and MT from ISA experiments to focus on the most tractable algorithms. 

\subsubsection{Details of Lagrangian Relaxation implementation}

Recall the derivation of LR for GAP in Section~\ref{sec:LR}. Reasonable effort has been made to ensure that our LR baseline is a fair comparison. We evaluated several approaches and selected the most efficient implementation; however, considering the vast literature on Lagrangian Relaxation, faster methods may exist.

At each iteration, we calculate the current lower bound objective and gradient by solving all knapsack subproblems.  For fair comparison, we use the same binary knapsack solver as our CG method. We terminate when the absolute change in objective is less than $10^{-6}$. The gradient is used to update the multipliers via a directional step; specifically, we use L-BFGS~\cite{LBFGS} with a memory of 256.  L-BFGS uses previously generated gradients to more accurately approximate the Hessian matrix for faster convergence.

\begin{enumerate}
    \item Initialize multipliers $\pi = \mathbf{0}$
    \item Solve knapsack subproblems for given $\pi$
    \item Compute lower bound objective and gradient from subproblem solution
    \item Update $\pi$ using L-BFGS backtracking line search step and gradient history
    \item Repeat until objective converges (or time out).
\end{enumerate}

\subsection{Termination and Integrality}

Our computational study terminates each instance at the CG root node; we do not solve the integer problem via branch-and-price. Introducing the complexity of branching would distract from our comparison on pricing strategies. Nevertheless, we can exploit the integrality of GAP to terminate early at the root node. Specifically, all Yagiura and ISA instances have integer cost coefficients, so the optimal objective must be integer and any valid lower bound can be tightened by rounding up. Formally, let the sum of negative reduced-costs be given by 
$$\text{rc}_t = \sum_{i \in I}\min\left\{rc\left(i\right), 0\right\}$$ 
where $rc(i)$ is the reduced cost for the $i^{\text{\,th}}$ machine as defined in \eqref{eq:dantzig-pricing}, and let $\text{RMP}_t$ be the objective value of the RMP in the $t^{\text{\,th}}$ iteration. 
Since $\text{rc}_t$ is an optimistic estimate of the remaining improvement in $\text{RMP}_t$ and we have convexity constraints \eqref{eq:mpcon2}, we have that $\text{RMP}_t + \text{rc}_t$ is a valid lower bound~\cite{lubbecke_column_2011}.
In our implementation, termination based on lower bounds occurs either when $|\text{rc}_t| < 10^{-6}$, or when $\lceil \max_t \{\text{RMP}_t + \text{rc}_t \} \rceil \ge \text{RMP}_t$ by exploiting integrality.  

Note that primal heuristics are not explicitly implemented, so a feasible integer upper bound is only obtained when the RMP solution happens to be integral at an iteration. Empirically, this rarely occurs with Dantzig or Pessoa, but is common with Template pricing. Regardless, consider $\text{UB}_t$, the best upper bound found at iteration $t$, we also terminate when the MIP gap is less than 0.001\%, that is: $\text{UB}_t - \lceil \max_t \{\text{RMP}_t + \text{rc}_t \} \rceil < 0.00001\, \text{UB}_t$.

\subsection{Initialization}


Since CG is inherently a primal method, it requires a primal feasible RMP at the start. Typical initialization approaches include heuristic, Phase I (artificial variables with big M), and Farkas pricing. These approaches are discussed by \cite{lubbecke_column_2011,vanderbeck_implementing_2005}. 
There are benefits and drawbacks to each; we briefly justify our choice.

Heuristics often yield high quality integer solutions; however, may struggle to find feasibility on challenging instances, and are often not ideal starting points for Dantzig or Pessoa pricing (due to integrality).
Phase I with big M, requires suitably tuning large costs per instance for best performance. 
Farkas pricing~\cite{achterberg2009scip,lubbecke_column_2011}, uses the Farkas certificate of infeasibility as the dual values. It avoids tuning, but the construction of dual values is dependent on the implementation underlying LP solver.

Our choice for initialization is equivalent to Farkas pricing, but is best described as Phase I (\textit{without} big M).
Specifically, we solve the following RMP; adding artificial variables $y^+_j$ and $y^-_j$ for $j \in J$ and minimizing their sum.  When the optimal objective of this modified RMP is non-zero, the dual variables $\pi$ and $\mu$ give the corresponding Farkas certificate for the original RMP. 
\begin{center}
    \renewcommand*{\arraystretch}{1.5}
    \begin{tabularx}{\textwidth}{>{\raggedleft\arraybackslash}rr@{\hskip 5pt}c@{\hskip 5pt}ll@{\hskip 5pt}l}
		$\displaystyle \min_{\lambda, y}$ & \multicolumn{4}{l}{
			$\;\displaystyle \sum_{j\in J} \left(y^+_j + y^-_j\right)$} & \\
		
		s.t. & $\;\displaystyle \sum_{i \in I} \sum_{\mathbf{v}_{ip} \in P_i} \mathbb{1}_{\{\mathbf{e}_j^\top\mathbf{v}_{ip} = 1\}} \lambda_{ip}$ & $\geq$ & $1  + \left(y^+_j - y^-_j\right)$ & $\quad\forall j\in J $ & $(\pi_j)$ \\
		
		& $\;\displaystyle \sum_{\mathbf{v}_{ip} \in P_i} \lambda_{ip}$ & $=$ & $1$ & $\quad\forall i\in I$ & $(\mu_i)$ \\
		
		& $\;\lambda_{ip}$ & $\geq$ & $ 0$ & $\quad\forall i\in I,\; \mathbf{v}_{ip}\in P_i$ & \\
		& $\; y^+_j, y^-_j$ & $\geq$ & $ 0$ & $\quad\forall j\in J.$ &
    \end{tabularx}
\end{center}
For the above, recall that in Section~\ref{sec:algorithms} our implementation uses a set cover formulation.

Since the modified RMP ignores column costs, the pricing subproblems also ignore the costs.  Specifically, Dantzig and Template pricing become:
\begin{equation*}
\min_{x\in P_i} \, -\pi^\intercal x, \;\text{and} \min_{x\in S^i_{(\pi,\mu)}} \, \left(-d\left(x,y^i\right),\; -\pi^\intercal x\right).
\end{equation*}
As far as we are aware, Phase I initialization (Farkas pricing) traditionally uses Dantzig pricing until reaching RMP feasibility. While it might be possible to adapt Pessoa pricing to be used during initialization, we do not attempt it in our study. For our Pessoa results, we first initialize the RMP using Dantzig pricing.
Template pricing can be used during Phase I initialization with a minor change. Templates can be easily constructed from the solution of the modified RMP; however, for the initial template (when the modified RMP is initially empty), our implementation uses the LP solution of the compact GAP formulation. Heuristics or other methods could also be used. 

The modified RMP is solved iteratively; generating columns every iteration until the optimal objective becomes zero -- indicating that a feasible solution to the original RMP exists. At this point, we can remove all artificial variables, and reinstate the correct costs for each column in the RMP.  Since the RMP is now initialized with a feasible solution, we continue to solve to optimality.

\begin{table}[hbtp]
\centering
\newcolumntype{d}{R{14mm}{4618}{}{3.5mm}}
\newcolumntype{i}{R{11mm}{6}{Dandelion}{2.5mm}}
\newcolumntype{g}{R{11mm}{201}{Peach}{2.5mm}}
\newcolumntype{e}{R{11mm}{3}{Red!65}{2.5mm}}

\begin{tabular}{ligedddd}
 & \multicolumn{3}{c}{Initialization} & \multicolumn{4}{c}{Convergence to optimality (\#its)} \\  \cmidrule(l){2-4} \cmidrule(l){5-8}
Phase I & \#its & \% gap & $\times$basis & Dantzig & Pessoa &    LT & MT \EndTableHeader \\ \hline\toprule
Dantzig    &  5 & 200.8 & 2.9 & 3469 &  634 & 2104 & 2315  \\
LT &  6 &   2.0 & 1.5 & 4439 & 1296 &  288 &  368  \\
MT &  5 &   2.0 & 1.0 & 4618 & 1314 &  295 &  370 
\end{tabular}

\vspace{1em}
\caption{Comparison of initialization quality and convergence on the c10400 instance (seed = 0). The $\times$basis indicates the fractionality of the initial solution (1.0 is integral). $\%$ gap stands for the quality of the initial feasible solution discovered compared to the optimal RMP objective value. $\# $its is the number of CG iterations. See Figure \ref{fig:c10400-convergence} for convergence behavior.}
\label{tab:c10400-farkas-its}
\end{table}

\begin{figure}[htbp]
\centering
\begin{tikzpicture}
\begin{axis}[
  name=left fig,
  axis x line*=bottom,
  xmin=0, xmax=4700,
  ymin=-2, ymax=210,
  ylabel={\% gap},
  ylabel near ticks,
  ylabel style={at={(axis description cs:-0.07,-0.2)}},
  tick label style={font=\small}, 
  width=0.95\columnwidth,
  height=5cm,
  legend cell align={left},
  legend style={/tikz/every even column/.append style={column sep=0.5cm}, font = \small, at={(0.5,1.1)},legend cell align={left}, anchor=north,legend columns=-1, draw=none, fill=none}
] 

\addplot [draw=mymt, line width=3.5pt] 
    table[x=iterations, y expr=100*\thisrow{d-MT}] {c10400-convergence.tsv};

\addplot [draw=mylt, line width=1.5pt] 
    table[x=iterations, y expr=100*\thisrow{d-LT}] {c10400-convergence.tsv};

\addplot [draw=myp, line width=2.0pt, cheating dash=on 4.0pt off 1.5\pgflinewidth, line cap=round, line join=round] 
    table[x=iterations, y expr=100*\thisrow{d-p}] {c10400-convergence.tsv};

\addplot [draw=myd, line width=3.5pt, cheating dash=on 0pt off 2.0\pgflinewidth, line cap=round, line join=round] 
    table[x=iterations, y expr=100*\thisrow{d-d}] {c10400-convergence.tsv};

\node[anchor=east] at ($ (axis cs:4500,-1) + (0cm,0.85cm) $) {\footnotesize {\textbf{Dantzig Phase I}}};

\legend{MT, LT, Pessoa, Dantzig}
\end{axis}

\begin{axis}[
  name=left fig,
  axis x line*=bottom,
  xmin=0, xmax=4700,
  ymin=-0.05, ymax=2.2,
  tick label style={font=\small}, 
  width=0.95\columnwidth,
  height=3cm,
  legend cell align={left},
  shift={(0, -2.1cm)},
] 
\addplot [draw=mymt, line width=3.5pt] 
    table[x=iterations, y expr=100*\thisrow{LT-MT}] {c10400-convergence.tsv};

\addplot [draw=mylt, line width=1.5pt] 
    table[x=iterations, y expr=100*\thisrow{LT-LT}] {c10400-convergence.tsv};

\addplot [draw=myp, line width=2.0pt, cheating dash=on 4.0pt off 1.5\pgflinewidth, line cap=round, line join=round] 
    table[x=iterations, y expr=100*\thisrow{LT-p}] {c10400-convergence.tsv};

\addplot [draw=myd, line width=3.5pt, cheating dash=on 0pt off 2\pgflinewidth, line cap=round, line join=round] 
    table[x=iterations, y expr=100*\thisrow{LT-d}] {c10400-convergence.tsv};
    
\node[anchor=east] at ($ (axis cs:4500,-0.05) + (0cm,0.68cm) $) {\footnotesize {\textbf{LT Phase I}}};

\end{axis}

\begin{axis}[
  name=left fig,
  axis x line*=bottom,
  xmin=0, xmax=4700,
  ymin=-0.05, ymax=2.2,
  xlabel={\# iterations},
  tick label style={font=\small}, 
  width=0.95\columnwidth,
  height=3cm,
  legend cell align={left},
  shift={(0, -4.1cm)},
] 
\addplot [draw=mymt, line width=3.5pt] 
    table[x=iterations, y expr=100*\thisrow{MT-MT}] {c10400-convergence.tsv};

\addplot [draw=mylt, line width=1.5pt] 
    table[x=iterations, y expr=100*\thisrow{MT-LT}] {c10400-convergence.tsv};

\addplot [draw=myp, line width=2.0pt, cheating dash=on 4.0pt off 1.5\pgflinewidth, line cap=round, line join=round] 
    table[x=iterations, y expr=100*\thisrow{MT-p}] {c10400-convergence.tsv};

\addplot [draw=myd, line width=3.5pt, cheating dash=on 0pt off 2\pgflinewidth, line cap=round, line join=round] 
    table[x=iterations, y expr=100*\thisrow{MT-d}] {c10400-convergence.tsv};
    
\node[anchor=east] at ($ (axis cs:4500,-0.05) + (0cm,0.68cm) $) {\footnotesize {\textbf{MT Phase I}}};

\end{axis}
\end{tikzpicture} 
\caption{Convergence of c10400 to optimality with various initialization methods (seed = 0).  See Table \ref{tab:c10400-farkas-its} for additional statistics.}
\label{fig:c10400-convergence}
\end{figure}
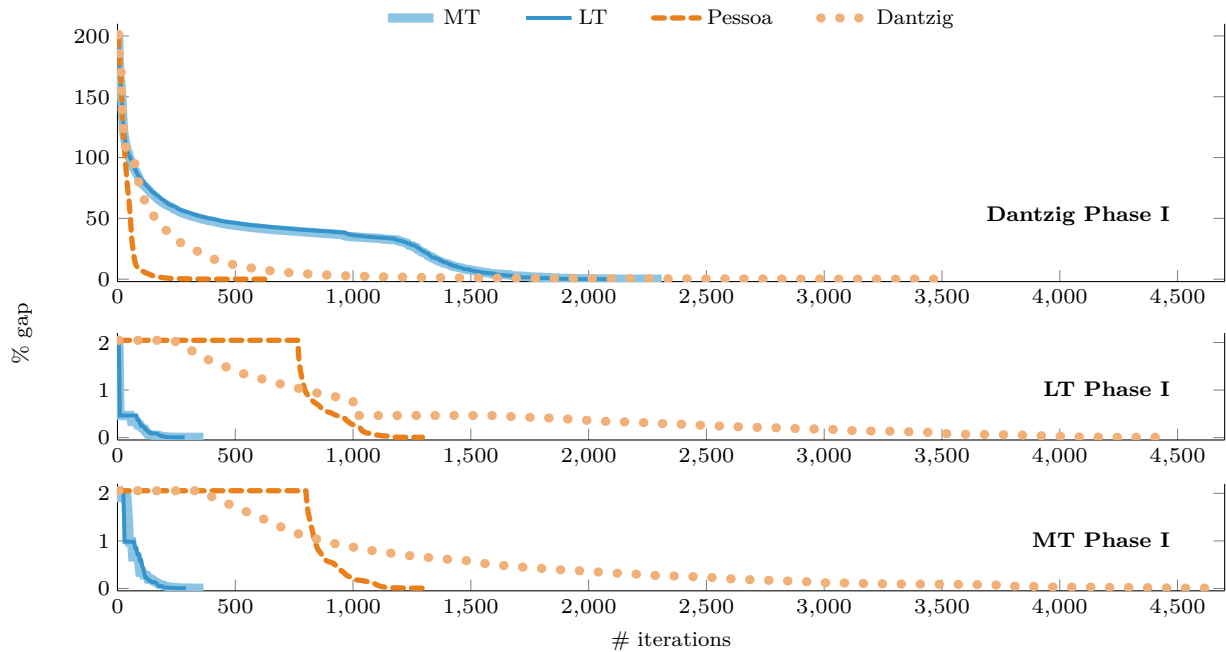

Initialization of the RMP has a significant impact on convergence. For our experiments to be as fair as possible, we consider this impact. Table~\ref{tab:c10400-farkas-its} shows (for a single instance) that Template Phase I finds a near optimal integer initial solution; whereas Dantzig Phase I has significantly worse quality and high fractionality. One might think that using Template Phase I to initialize all methods gives a fair comparison; Table \ref{tab:c10400-farkas-its} and Figure \ref{fig:c10400-convergence} show that this is not the case.  In particular, the integrality of Template Phase I seems to slow convergence for Dantzig and Pessoa pricing methods to converge to optimality; whereas the fractionality of Dantzig Phase I hurts Template pricing methods convergence to the optimal RMP solution.  Using integer solutions as a starting basis has long been known~\cite{vanderbeck_decomposition_1994} to yield poor convergence, so this behavior of Template pricing is interesting.

Note that when the RMP is highly fractional, our chosen template scheme may result in many or all elements set to zero. When this occurs, Template pricing essentially reverts to Dantzig pricing; however, our dynamic scheme implicitly tries to recover the template structure and terminates sooner than Dantzig.  Other template schemes could alternatively revert to Pessoa instead of Dantzig; we have chosen this particular scheme for its simplicity.

\subsection{Column management}
\label{sec:column-management}

Our implementation aims to add one column per machine at each iteration. After many iterations, however, the RMP can become unnecessarily large. For example, with the instance in Table \ref{tab:c10400-farkas-its}, after 4000 iterations, the RMP could contain 40,000 columns; whereas the RMP basis requires at most 410 columns, and an integer solution could be represented using only 10 columns.

RMP solve times increase with the number of columns, many of which are not needed for convergence. Thus, we periodically remove columns to improve performance. Identifying the minimal set of necessary columns is difficult, so it is performed heuristically. Common choices include age and reduced-cost based retention policies; we use age-based retention for its simplicity.

Age-based retention keeps, for each column, the most recent iteration in which it was in the basis. At each iteration, only basic columns need to be updated. Newly added columns are initialized with the iteration in which they were added, regardless of whether they immediately enter the basis. Retention is controlled by the age threshold $\tau$. At iteration $t$, we retain columns whose recorded age is in the interval $\{t-\tau, \ldots, t\}$, and remove all others.

Larger values of $\tau$ retain more columns in the RMP, making it less likely that columns important for approximating the dual polyhedron are removed. Small values of $\tau$ can occasionally produce very fast solves when the retained columns happen, by chance, to be sufficient; however, this behavior is unreliable and typically leads to highly volatile performance due both to dual instability and to the repeated regeneration of previously discarded columns. Increasing $\tau$ reduces these destabilizing effects and yields more robust performance, at the cost of larger RMPs and increased solve times. We seek a policy that minimizes the RMP solve time while remaining robust, in the sense that small changes in the threshold do not significantly impact performance.

\begin{figure}[hbtp]
\centering
\begin{tikzpicture}
\definecolor{myd}{HTML}{f0b077}
\definecolor{myp}{HTML}{ea801c}
\definecolor{mylt}{HTML}{3594cc}
\definecolor{mymt}{HTML}{8cc5e3}

\tikzset{
  eqnlabel/.style={
    inner sep=1pt,
    text opacity=1,
    font=\tiny
  }
}

\begin{axis}[
  name=left fig,
  axis x line*=bottom,
  xmin=0, xmax=85,
  ymin=0, ymax=140,
  ylabel={Best age threshold},
  ylabel near ticks,
  xlabel={Job to Machine Ratio},
  tick label style={font=\small}, 
  width=0.95\columnwidth,
  height=4.5cm,
  legend cell align={left},
  legend style={/tikz/every even column/.append style={column sep=0.5cm}, font = \small, at={(0.5,1.3)},legend cell align={left}, anchor=north,legend columns=-1, draw=none, fill=none}
]
    \addplot [mylt, line width=1.0pt, domain=0:85, samples=201] (
        x,{0.00044*x^2 + 0.0405*x + 1}
    );
    \addplot [myp, line width=2.0pt, cheating dash=on 4.0pt off 1.5\pgflinewidth, line cap=round, line join=round, domain=0:85, samples=201] (
        x,{0.0*x^2 + 0.3*x + 1}
    );
    \addplot [myd, line width=3.5pt, cheating dash=on 0pt off 2\pgflinewidth, line cap=round, line join=round, domain=0:85, samples=201] (
        x,{0.081875*x^2 - 0.0*x + 1}
    );

\addplot [draw=myd, only marks, mark=*, mark options={fill=myd,solid,fill opacity=0.1,draw opacity=1, scale=0.6}] 
    table[x=job_machine_ratio, y=best_age] {age-sweep-dantzig.tsv};

\addplot [draw=myp, only marks, mark=x, mark options={scale=1}] 
    table[x=job_machine_ratio, y=best_age] {age-sweep-pessoa.tsv};

\addplot [draw=mylt, only marks, mark=*, mark options={fill=white, scale=0.7}] 
    table[x=job_machine_ratio, y=best_age] {age-sweep-lt.tsv};

\addplot[draw=none, domain=0:85, samples=201] {0.081875*x^2 + 1}
  node[pos=0.1, sloped, transform shape, above, yshift=3pt, eqnlabel] {$0.081875x^2+1$};

\addplot[draw=none, domain=0:85, samples=201] {0.3*x + 1}
  node[pos=0.55, sloped, transform shape, above, yshift=3pt, eqnlabel] {$0.3x + 1$};

\addplot[draw=none, domain=0:85, samples=201] {0.00044*x^2 + 0.0405*x + 1}
  node[pos=0.87, sloped, transform shape, above, yshift=1pt, eqnlabel] {$0.00044x^2 + 0.0405x + 1$};
    
\legend{LT, Pessoa, Dantzig}
\end{axis}
\end{tikzpicture} 
\caption{Chosen age threshold policy vs. degeneracy and algorithm. Policy curves are determined by minimizing the quadratic that covers all data points.}
\label{fig:retention-sweep}
\end{figure}
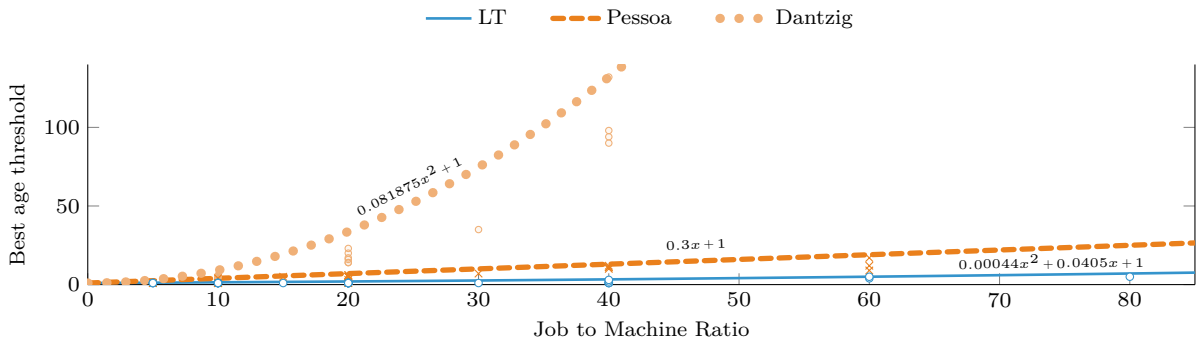

We also considered an activation-frequency parameter controlling when column removal is triggered; however, preliminary experiments indicated that activating at every iteration was sufficient in our setting, with little downside under primal simplex since the basis remains unchanged. We therefore fix activation at every iteration and determine the age threshold for each instance using the policy shown in Figure~\ref{fig:retention-sweep}. The chosen threshold reflects an empirically selected near-optimal choice for each pricing rule and job to machine ratio combination in order to reduce the time for CG convergence. Details on how these policies are determined via an extensive parameter sweep are presented in Appendix \ref{app:column-management}.

An important takeaway is that stronger dual stabilization reduces sensitivity to age threshold. In particular, Template pricing appears to require the least tuning and is the most robust among the methods we tested.

\section{Experimental Results} \label{sec:results}

We present computational experiments that evaluate Template pricing and compare with other approaches.
All methods are implemented in C++, with HiGHS 1.11.0 \cite{huangfu_parallelizing_2018} used to solve the RMP. All pricing algorithms share the same knapsack solver to ensure fair comparison -- except MT, which solves an IP. The knapsack solver implements a standard dynamic programming algorithm, with the code adapted from SCIP \cite{SCIPOptSuite10}.
The experiments were run on an AMD EPYC 7V12 64-Core 1.5 GHz processor with 1TB RAM. Each instance had exclusive access to the hardware and the subproblems were run in parallel.

The details on instances, initialization, termination criteria, and column management have been described in the previous section. Further details of the implementation can be found in the released source code \cite{templatecode}.

\subsection{Yagiura Results}

The experiments on Yagiura use a time limit of 6 hours and have 5 replications (different random seed) per instance, unless otherwise specified.  A total of 285 experiments per algorithm. We evaluate Dantzig (D), Pessoa (P), Lagrange Template (LT), and MIP Template (MT) pricing.

\subsubsection{Yagiura Initialization}

Recall that our RMP uses a set cover instead of a set partition formulation for additional stabilization; moreover, with this choice of formulation, it is also much easier to construct an initial feasible RMP -- although, by design, this may include jobs more than once. We compare the impact of formulation choice on Phase I initialization in Table \ref{tab:yagiura-init}.

\begin{table}[htbp]
\centering
\newcolumntype{t}{R{18mm}{172988}{Dandelion}{1.5mm}}
\newcolumntype{i}{R{10mm}{77}{Peach}{1.5mm}}
\newcolumntype{g}{R{11mm}{337}{Peach}{1.5mm}}
\newcolumntype{e}{R{15mm}{100}{YellowGreen}{1.5mm}}

\begin{subtable}[t]{\textwidth}
\centering

\begin{tabular}{lttigee}
Phase I & RMP (ms) & Pricing (ms) & \# its & \% gap & \% integral & \% timeout \EndTableHeader \\ \hline\toprule
Dantzig  & 18.1 &   2.1 & 4.4 & 336.9 & 31.9 & 0.0 \\
LT       &  2.3 &   4.2 & 3.8 &   6.5 & 83.5 & 0.0 \\
MT       &  4.0 & 117.6 & 4.8 &   5.0 & 99.3 & 0.0
\end{tabular}
\vspace{1em}
\caption{RMP set cover}
\end{subtable}

	\vspace{0.75em}

\begin{subtable}[t]{\textwidth}
\centering

\begin{tabular}{lttigee}
Phase I & RMP (ms) & Pricing (ms) & \# its & \% gap & \% integral & \% timeout \EndTableHeader \\ \hline\toprule
Dantzig  & 172987.9 &  69.2 &  76.2 & 103.7 &  0.0 & 26.3 \\
LT       &  34936.6 & 128.4 &  33.8 &  16.3 & 49.5 &  3.2 \\
MT       &     56.2 & 603.6 &  16.9 &   6.1 & 97.5 &  0.0 
\end{tabular}
\vspace{1em}
\caption{RMP set partition}
\end{subtable}

\caption{Yagiura: Comparison of initialization quality. \% gap is relative to optimal RMP (excludes timeout). 5 replications, 10 min time limit.}
\label{tab:yagiura-init}
\end{table}


For these experiments we have restricted the time limit to 10 minutes. 
As can be seen by these results, Template Phase I (both LT and MT) finds high quality initial solutions in similar or fewer iterations; while taking less compute time in the RMP. The Template solutions are also more likely to be integral. The differences between LT and MT are emphasized by the set partition results; the benefits of exact Template Phase I (MT) can be clearly seen, however, heuristic Template Phase I (LT) still dominates over Dantzig Phase I.

\subsubsection{Yagiura Results}\label{sec:yagiura-results}

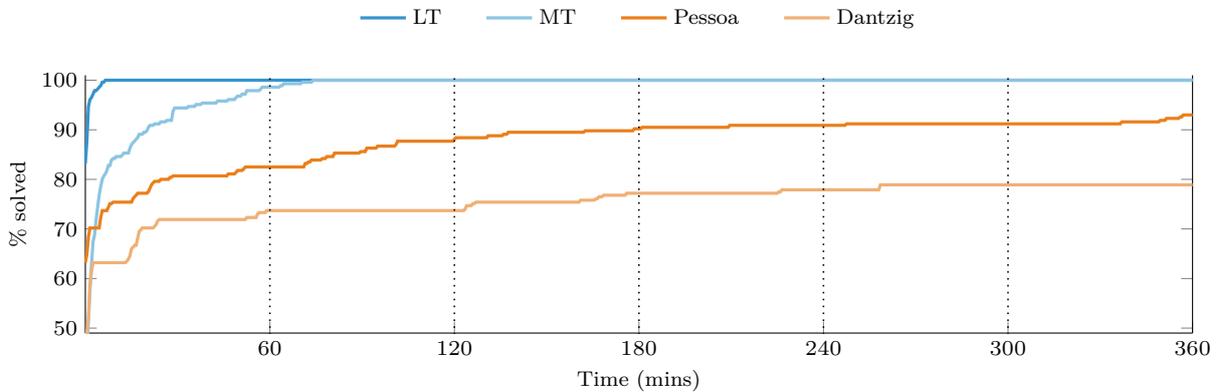
\begin{figure}[htbp]
\centering
\definecolor{myd}{HTML}{f0b077}
\definecolor{myp}{HTML}{ea801c}
\definecolor{mylt}{HTML}{3594cc}
\definecolor{mymt}{HTML}{8cc5e3}

\begin{tikzpicture}
\begin{axis}[
  name=left fig,
  axis x line*=bottom,
  xmin=0, xmax=360,
  ymin=49, ymax=101,
  ylabel={\% solved},
  ylabel near ticks,
  xlabel={Time (mins)},
  xtick={60,120,180,240,300,360},
  ytick={50, 60, 70, 80, 90, 100},
  tick label style={font=\small}, 
  width=0.95\columnwidth,
  height=5cm,
  legend cell align={left},
  legend style={/tikz/every even column/.append style={column sep=0.5cm}, font = \small, at={(0.5,1.3)},legend cell align={left}, anchor=north,legend columns=-1, draw=none, fill=none}
] 
\addplot [draw=mylt, very thick] 
    table[x=time, y=LT] {count-vs-time-yagiura.tsv};

\addplot [draw=mymt, very thick] 
    table[x=time, y=MT] {count-vs-time-yagiura.tsv};

\addplot [draw=myp, very thick] 
    table[x=time, y=P] {count-vs-time-yagiura.tsv};

\addplot [draw=myd, very thick] 
    table[x=time, y=D] {count-vs-time-yagiura.tsv};

\draw [dotted, semithick] (axis cs: 60, 0) -- (axis cs: 60, 103);
\draw [dotted, semithick] (axis cs: 120, 0) -- (axis cs: 120, 103);
\draw [dotted, semithick] (axis cs: 180, 0) -- (axis cs: 180, 103);
\draw [dotted, semithick] (axis cs: 240, 0) -- (axis cs: 240, 103);
\draw [dotted, semithick] (axis cs: 300, 0) -- (axis cs: 300, 103);
    
\legend{LT, MT, Pessoa, Dantzig}
\end{axis}
\end{tikzpicture} 
\caption{Yagiura: solved vs solve time.}
\label{fig:yagiura-count-vs-time}
\end{figure}

Figure \ref{fig:yagiura-count-vs-time} presents the percentage of instances that solved over the 6-hour time limit. Both LT and MT successfully solve all instances, clearly outperforming D and P pricing methods. In particular, LT solved each instance in less than 7 minutes. Summary statistics are provided in Table \ref{tab:yagiura-combined}. Notably, LT and MT requires approximately one-third as many iterations while consuming orders of magnitude less RMP time. This indicates that the reduction in RMP time is not solely attributable to fewer iterations. In particular, consider the number of pivots required per column and pivots per second; Template pricing requires around $10\times$ fewer pivots, each executed $10\times$ more efficiently. 

\begin{table}[htbp]
\centering
  
\begin{subtable}[t]{\textwidth}
	\centering
	\caption{Computational effort of RMP and pricing}

	\begin{tabular}{lR{25mm}{21.1}{YellowOrange!60}{2.5mm}R{20mm}{940}{YellowOrange!60}{2.5mm}R{20mm}{5290}{YellowOrange!60}{2.5mm}R{20mm}{5290}{YellowOrange!60}{2.5mm}}
	Method & \% timeout & \# CG its & RMP (s) & Pricing (s) \EndTableHeader \\ \toprule
	D  & 21.1 & 940 & 5285.2 &   1.0 \\
	P  &  7.0 & 438 & 2642.3 &   1.6 \\
	LT &  0.0 & 325 &   15.4 &   5.7 \\
	MT &  0.0 & 436 &   14.9 & 372.8 \\
	\end{tabular}
\end{subtable}

	\vspace{0.75em}

\begin{subtable}[t]{\textwidth}
	\centering
	\caption{Pivot efficiency and integrality of solutions}
	\begin{tabular}{lR{25mm}{53.5}{YellowOrange!60}{2.5mm}R{20mm}{2364}{ForestGreen!60}{2.5mm}R{18mm}{100}{ForestGreen!60}{2.5mm}R{20mm}{51}{YellowOrange!60}{2.5mm}}
	Method & RMP pivots/col & RMP pivots/s & \% integral & integer gap \% \EndTableHeader\\ \toprule
	D & 53.8 &  215.4 &  37.9 & 50.7 \\
	P & 37.9 &  222.1 &  40.0 & 46.7 \\
	LT & 3.6 & 2363.4 & 100.0 &  0.8 \\
	MT & 2.7 & 2019.3 & 100.0 &  0.3 \\
	\end{tabular}
\end{subtable}

\caption{Yagiura: high-level summary.}
\label{tab:yagiura-combined}
\end{table}

A possible explanation for fewer pivots and these pivots being faster with Template pricing is the following.
Template pricing aims to generate nearly-compatible sets of columns that are structurally similar to those in the current basis $B$. 
Consider the newly generated column $a_i$ for machine $i\in I$. When $a_i$ is structurally similar to the basic column $a_{B_r}$ we have $B^{-1}a_i \approx e_r$, that is, the unit vector with 1 in position $r$ (corresponding to that basic column).  The simplex algorithm calculates $B^{-1}a_i$ via FTRAN (forward transformation) and this linear algebra is fastest when both $a_i$ and the resulting vector is sparse; moreover, sparsity implies fewer candidates to evaluate in the ratio-test and ultimately yields a more localized basis change; that is, the first pivot is almost a one-for-one replacement.  After this pivot, $B^{-1}$ incurs only a small change, which has multiple benefits. First, this implies only a relatively small perturbation in the duals $c^{\intercal}_BB^{-1}$, which is less likely to trigger a cascade of additional improving pivots. At the same time, relatively small changes to $B^{-1}$ requires less bookkeeping for steepest-edge or Devex pricing-weight updates and fewer re-inversions or refactorization of $B$. So, generated columns that are structurally close to the current basis tend to create sparser transformed columns and sparser update matrices, which means each solve is cheaper, the accumulated basis-factor updates grow more slowly, and costly reinversion and refactorization can be postponed or made cheaper when it does occur.  Moreover, at each iteration we add a set of nearly-compatible columns, so the impact is even stronger; each pivot tends to replace different basic columns with limited interaction across rows. These new columns can enter the basis with relatively little ``undoing'' of previous pivots. From a linear algebra perspective, the pivotal columns remain sparse and separated, so FTRANs, ratio tests, pricing-weight updates, and basis-factor updates all stay relatively cheap.

Furthermore, Template pricing consistently finds high-quality integer solutions for all instances. Finally, Table \ref{tab:yagiura-rmp_cg_times} provides a breakdown of RMP and CG pricing time by job and machine, highlighting the effects of scale and degeneracy (approximated by the jobs/machines ratio).

\begin{table}[htbp]
\centering
\newcolumntype{d}{R{14mm}{21600}{}{2.5mm}}
\newcolumntype{g}{R{14mm}{4000}{}{2.5mm}}
\resizebox{0.7\columnwidth}{!}{%
\begin{tabular}{@{}rrddddgggg}
&& \multicolumn{4}{c}{RMP (s)} & \multicolumn{4}{c}{Pricing (s)} \\
\cmidrule(l){3-6}\cmidrule(l){7-10}
 $|J|$& $|I|$& D & P & LT & MT & D & P & LT & MT \EndTableHeader\\
\midrule
100 &  5 & 1.6 & 0.3 & 0.1 & 0.2 & 0.2 & 0.1 & 0.1 & 37.0 \\
    & 10 & 0.4 & 0.2 & 0.0 & 0.1 & 0.1 & 0.0 & 0.0 & 11.6 \\
    & 20 & 0.2 & 0.1 & 0.0 & 0.1 & 0.0 & 0.0 & 0.0 &  6.8 \\\midrule
     
200 &  5 & 96.5 & 6.7 & 0.8 & 1.0 & 1.3 & 0.5 & 0.9 & 108.6 \\
    & 10 & 12.5 & 1.8 & 0.3 & 0.4 & 0.2 & 0.2 & 0.4 &  45.6 \\
    & 20 &  3.0 & 1.0 & 0.2 & 0.3 & 0.1 & 0.1 & 0.1 &  25.5 \\\midrule
     
400 & 10 & 957.1 & 67.0 & 3.8 & 4.1 & 3.6 & 1.6 & 3.7 & 223.0 \\
    & 20 & 131.5 & 19.3 & 1.6 & 2.1 & 0.6 & 0.4 & 1.4 & 198.3 \\
    & 40 &  32.3 &  8.8 & 0.9 & 1.3 & 0.2 & 0.2 & 0.4 &  88.7 \\\midrule
 
900 & 15 & 21596.0 & 4147.4 & 26.3 & 42.2 & 3.0 & 11.7 & 16.2 & 1180.3 \\
    & 30 & 10769.0 &  864.7 &  9.8 & 12.8 & 4.6 &  3.0 &  7.0 &  600.1 \\
    & 60 &  1946.5 &  302.7 &  5.5 &  7.2 & 1.0 &  0.9 &  2.2 &  454.1 \\\midrule

1600 & 20 & 21600.0 & 21597.5 & 163.9 & 144.6 & 1.4 & 3.5 & 48.1 & 1398.7 \\
     & 40 & 21598.0 & 16788.3 &  46.7 &  42.9 & 1.4 & 5.3 & 15.7 & 1238.9 \\
     & 80 & 21598.0 &  6390.4 &  32.4 &  21.9 & 0.6 & 3.1 & 10.2 & 1308.8 \\     
\end{tabular}
}
\caption{Yagiura: RMP and Pricing time (s) for different job/machine}
\label{tab:yagiura-rmp_cg_times}
\end{table}

\subsection{ISA Instances}

The ISA experiments use a time limit of 1 hour and perform a single replication per instance. Thus, a total of 1735 experiments per algorithm were conducted. We primarily compare Pessoa (P) and Lagrange Template (LT), however also include results from MIP Template (MT) and Lagrangian Relaxation (LR) where appropriate.

\subsubsection{ISA Initialization}

We perform additional experiments to highlight the quality of initialization. Table \ref{tab:isa-init} shows a comparison of the various Phase I initialization with a time limit of 10 minutes.  Dantzig Phase I fails to find a primal feasible RMP within this time limit for 2.1\% of the instances. In a separate experiment, we found that some of these instances required over 74 hours to initialize a feasible RMP.  Template Phase I initialized all instances within the time limit; spending a maximum of around 1 minute for LT and 5 minutes for MT.

As expected, Template Phase I finds initial solutions with significantly lower gap, and are much more likely to be integral -- requiring less time and similar (or fewer) iterations.  The quality difference can be seen between the heuristic LT and the exact method MT.

\begin{table}[htbp]
\centering
\newcolumntype{t}{R{18mm}{21}{Dandelion}{1.5mm}}
\newcolumntype{i}{R{10mm}{8}{Peach}{1.5mm}}
\newcolumntype{g}{R{11mm}{1160}{Peach}{1.5mm}}
\newcolumntype{e}{R{15mm}{100}{YellowGreen}{1.5mm}}

\begin{tabular}{lttigee}
Phase I & RMP (s) & Pricing (s) & \# its & \% gap & \% integral & \% timeout \EndTableHeader \\ \hline\toprule
Dantzig  & 20.3 &  0.1 & 6.2 & 1156.5 & 34.8 &  2.1 \\
LT       &  0.1 &  0.1 & 4.2 &   17.3 & 68.2 &  0.0 \\
MT       &  0.8 &  1.7 & 7.9 &   11.2 & 89.0 &  0.0
\end{tabular}
\vspace{1em}
\caption{ISA instances: Comparison of initialization statistics on ISA instances. \% gap is relative to optimal RMP bound. 10 min time limit.}
\label{tab:isa-init}
\end{table}

\subsubsection{ISA Results}

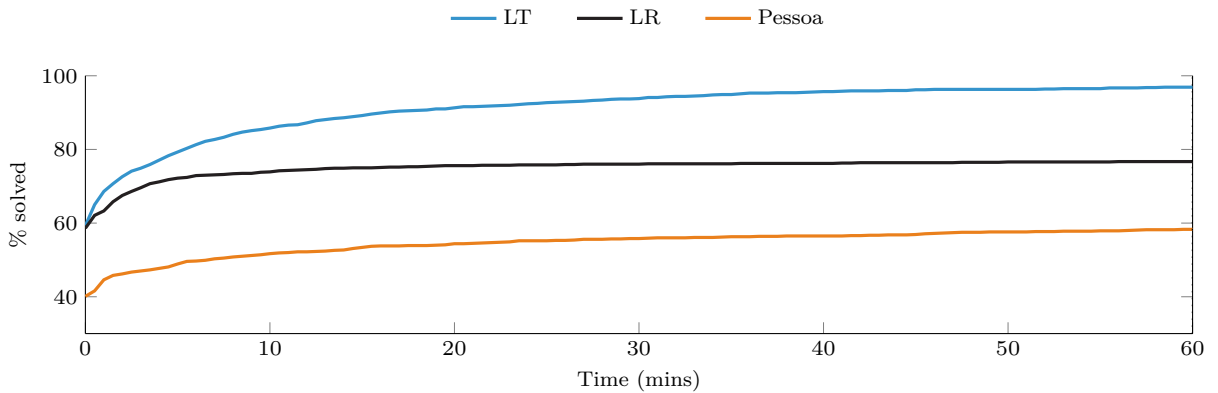
\begin{figure}[htbp]
\centering
\begin{tikzpicture}
\begin{axis}[
  name=left fig,
  axis x line*=bottom,
  xmin=0, xmax=60,
  ymin=30, ymax=100,
  ylabel={\% solved},
  ylabel near ticks,
  xlabel={Time (mins)},
  xtick={0,10,20,30,40,50,60},
  tick label style={font=\small}, 
  width=0.95\columnwidth,
  height=5cm,
  legend cell align={left},
  legend style={/tikz/every even column/.append style={column sep=0.5cm}, font = \small, at={(0.5,1.3)},legend cell align={left}, anchor=north,legend columns=-1, draw=none, fill=none}
] 
\addplot [draw=mylt, very thick] 
    table[x expr=\thisrow{time}, y expr=\thisrow{LT}] {count-vs-time-isa.tsv};

\addplot [draw=Black, very thick] 
    table[x expr=\thisrow{time}, y expr=\thisrow{LR}] {count-vs-time-isa.tsv};

\addplot [draw=myp, very thick] 
    table[x expr=\thisrow{time}, y expr=\thisrow{P}] {count-vs-time-isa.tsv};

\draw [dotted, semithick] (axis cs: 60, 0) -- (axis cs: 60, 103);
\draw [dotted, semithick] (axis cs: 120, 0) -- (axis cs: 120, 103);
\draw [dotted, semithick] (axis cs: 180, 0) -- (axis cs: 180, 103);
\draw [dotted, semithick] (axis cs: 240, 0) -- (axis cs: 240, 103);
\draw [dotted, semithick] (axis cs: 300, 0) -- (axis cs: 300, 103);
    
\legend{LT, LR, Pessoa}
\end{axis}
\end{tikzpicture} 
\caption{ISA: Solved vs solve time. Instances where LR does not reach optimal bound are treated as timed-out.}
\label{fig:isa-count-vs-time}
\end{figure}

Figure \ref{fig:isa-count-vs-time} shows the percentage of instances that solved over the 1-hour time limit.  Many of these instances are larger than those from Yagiura, and some were intentionally constructed to be more challenging by making feasible integer solutions harder to find.  LT performs well, solving nearly all instances; whereas LR quickly finds good bounds but fails to converge to the optimal. Pessoa solves 58\% of the instances. Extrapolating from the results in Yagiura, we expect Dantzig would perform considerably worse. See Table \ref{tab:isa-combined} for detailed results. In particular, we note that the RMP pivots are faster and fewer for LT than that of Pessoa, as also observed and discussed for the Yagiura instances in Section~\ref{sec:yagiura-results}.  Similarly, Table \ref{tab:isa-combined-solved} shows the results only on the 824 instances that were solved by all methods.  These instances can be considered `easy', and we expect the impact of LT is more pronounced in the harder instances.

In Table \ref{tab:isa-combined}, LR does not actually have a RMP nor generate columns per se, but there is a similar analogue, since at each iteration it solves the pricing problem for each machine (Pricing), and uses this information to update the multipliers (RMP). LR spends most of its time in CG pricing, with significantly more iterations; note these iterations (\# CG its) include the backtracking linesearch within L-BFGS. 
Although LR uses the same pricing algorithm, 0.5\% of the instances could not be solved by LR due to memory errors in the DP subproblem; these issues were not seen in the CG methods. 

\begin{table}[htbp]
\centering
  
\begin{subtable}[t]{\textwidth}
	\centering
	\caption{Computational effort of RMP and pricing}

	\begin{tabular}{lR{25mm}{100}{ForestGreen!60}{2.5mm}R{20mm}{7768}{YellowOrange!60}{2.5mm}R{20mm}{1631}{YellowOrange!60}{2.5mm}R{20mm}{1631}{YellowOrange!60}{2.5mm}}
	Method & \% solved & \# CG its & RMP (s) & Pricing (s) \EndTableHeader \\ \toprule
	P  &  58.3 &  245 & 1631.0 &   3.1 \\
	LT &  96.9 &  422 &  246.3 &  81.0 \\
    LR &  76.7 & 7768 &    2.9 & 409.0
	\end{tabular}
\end{subtable}

	\vspace{0.75em}

\begin{subtable}[t]{\textwidth}
	\centering
	\caption{Pivot efficiency and integrality of solutions}
	\begin{tabular}{lR{25mm}{54}{YellowOrange!60}{2.5mm}R{20mm}{956}{ForestGreen!60}{2.5mm}R{18mm}{100}{ForestGreen!60}{2.5mm}R{20mm}{42}{YellowOrange!60}{2.5mm}}
	Method & RMP pivots/col & RMP pivots/s & \% integral & integer gap \% \EndTableHeader\\ \toprule
	P & 53.3 & 310.0 & 39.6 & 41.3 \\
	LT & 5.1 & 955.7 & 98.3 &  1.3 \\
    LR &   - &     - &  6.5 &  0.8
	\end{tabular}
\end{subtable}

\caption{ISA: high-level statistics for all 1735 instances}
\label{tab:isa-combined}
\end{table}

\begin{table}[htbp]
\centering
  
\begin{subtable}[t]{\textwidth}
	\centering
	\caption{Computational effort of RMP and pricing}

	\begin{tabular}{lR{25mm}{100}{ForestGreen!60}{2.5mm}R{20mm}{1531}{YellowOrange!60}{2.5mm}R{20mm}{231}{YellowOrange!60}{2.5mm}R{20mm}{231}{YellowOrange!60}{2.5mm}}
	Method & \% solved & \# CG its & RMP (s) & Pricing (s) \EndTableHeader \\ \toprule
	P  &  100.0 &  341 & 230.9 &  1.1 \\
	LT &  100.0 &  226 &   3.0 &  2.3 \\
    LR &  100.0 & 1531 &   0.2 &  4.4
	\end{tabular}
\end{subtable}

	\vspace{0.75em}

\begin{subtable}[t]{\textwidth}
	\centering
	\caption{Pivot efficiency and integrality of solutions}
	\begin{tabular}{lR{25mm}{23}{YellowOrange!60}{2.5mm}R{20mm}{7075}{ForestGreen!60}{2.5mm}R{18mm}{100}{ForestGreen!60}{2.5mm}R{20mm}{30}{YellowOrange!60}{2.5mm}}
	Method & RMP pivots/col & RMP pivots/s & \% integral & integer gap \% \EndTableHeader\\ \toprule
	P  & 23.0 & 1227.4 & 45.6 & 29.8 \\
	LT &  2.9 & 7074.8 & 99.9 &  0.9 \\
    LR &    - &     - &  12.4 &  0.9
	\end{tabular}
\end{subtable}

\caption{ISA: high-level statistics for the 824 `easier' instances that were solved by all algorithms}
\label{tab:isa-combined-solved}
\end{table}

LR seldom found integer solutions implicitly, whereas LT consistently produces high-quality solutions across almost all instances.  For the  7\% of instances where LR obtained integer solutions, their quality was high; however, these instances were trivial to solve.  This can be seen in Table \ref{tab:isa-rmp_cg_times}, which explores how class and scale (number of jobs) impact outcomes.  Large number of jobs correlate with higher degeneracy, although difficulty also varies between classes. The `f' class instances, specifically designed to impede finding integer feasible solutions, proved challenging for all algorithms, though LT managed them best. As mentioned, LR only found integral solutions in trivial instances, that is, the `easier' classes with smaller number of jobs. As a dual method, LR generally does not produce feasible solutions without additional heuristics. In Pessoa’s case, achieving integrality appears random and unrelated to scale or degeneracy; however it failed to find any feasible integer solutions for the f class instances, and surprisingly also struggled on c class.

When LR does not converge to the optimal bound, we treat it as a timeout. Table~\ref{tab:isa-rmp_cg_times} highlights when this is the case, as it shows many timeouts even though the time to solve is relatively small. This accounts for 15\% of instances, where 13\% converged to a near optimal bound, while the remaining were far from optimal (the n class in particular).  LR is known to quickly give excellent bounds, while having issues with convergence; this is again shown here. Even considering this, LT still dominates with its stronger convergence guarantees, but it can take longer to reach similar bounds. A hybrid LR/LT approach could be very promising.

\begin{table}[htbp]
\centering
\resizebox{0.85\columnwidth}{!}{%
\newcolumntype{d}{R{12mm}{3600}{}{2.0mm}}
\newcolumntype{g}{R{12mm}{3600}{}{2.0mm}}
\newcolumntype{i}{R{8mm}{100}{ForestGreen!40}{2.5mm}}
\newcolumntype{t}{R{8mm}{100}{YellowOrange!60}{2.5mm}}

\begin{tabular}{@{}l@{}rdddgggiiittt}
&& \multicolumn{3}{c}{RMP (s)} & \multicolumn{3}{c}{Pricing (s)} & \multicolumn{3}{c}{\% integral} & \multicolumn{3}{c}{\% timeout}\\
\cmidrule(l){3-5}\cmidrule(l){6-8}\cmidrule(l){9-11}\cmidrule(l){12-14}
 class& $|J|$& P & LT & LR & P & LT & LR& P & LT & LR & P & LT & LR \EndTableHeader\\
\midrule
a & 100 & 0.2 & 0.0 & 0.0 & 0.0 & 0.0 & 0.0 & 87 & 100 & 43 & 0 & 0 & 0 \\
  & 200 & 2.2 & 0.0 & 0.0 & 0.1 & 0.0 & 0.0 & 63 & 100 & 13 & 0 & 0 & 0 \\\midrule
b & 100 & 0.2 & 0.0 & 0.0 & 0.1 & 0.0 & 0.0 & 57 & 100 & 30 & 0 & 0 & 10 \\
  & 200 & 2.8 & 0.3 & 0.0 & 0.2 & 0.2 & 0.1 & 20 & 100 & 10 & 0 & 0 & 10 \\\midrule

c & 100 & 0.2 & 0.0 & 0.0 & 0.0 & 0.0 & 0.1 & 13 & 100 & 17 & 0 & 0 & 13 \\
  & 200 & 2.7 & 0.3 & 0.0 & 0.2 & 0.2 & 0.1 & 7 & 100 & 7 & 0 & 0 & 10 \\
  & 400 & 18.9 & 0.7 & 0.1 & 0.3 & 0.3 & 0.2 & 0 & 100 & 0 & 0 & 0 & 3 \\
  & 900 & 1048.0 & 7.3 & 0.1 & 1.2 & 1.5 & 0.4 & 0 & 100 & 0 & 0 & 0 & 3 \\
  & 1600 & 2903.6 & 42.0 & 0.4 & 0.6 & 3.8 & 1.2 & 0 & 100 & 0 & 67 & 0 & 0 \\\midrule

d & 100 & 0.2 & 0.1 & 0.0 & 0.1 & 0.2 & 0.2 & 43 & 100 & 0 & 0 & 0 & 10 \\
 & 200 & 4.0 & 1.1 & 0.0 & 0.6 & 1.6 & 0.4 & 30 & 100 & 0 & 0 & 0 & 3 \\
 & 400 & 39.0 & 3.6 & 0.1 & 1.5 & 4.1 & 1.6 & 30 & 100 & 0 & 0 & 0 & 3 \\
 & 900 & 1783.4 & 19.4 & 0.4 & 5.3 & 15.0 & 14.7 & 47 & 100 & 0 & 33 & 0 & 0 \\
 & 1600 & 3597.6 & 114.6 & 1.0 & 2.6 & 74.0 & 73.4 & 47 & 100 & 0 & 100 & 0 & 3 \\\midrule

e & 100 & 0.2 & 0.1 & 0.0 & 0.1 & 0.0 & 0.1 & 43 & 100 & 0 & 0 & 0 & 37 \\
 & 200 & 3.9 & 0.5 & 0.0 & 0.2 & 0.3 & 0.2 & 37 & 100 & 7 & 0 & 0 & 13 \\
 & 400 & 36.6 & 1.8 & 0.1 & 0.4 & 0.9 & 0.3 & 50 & 100 & 0 & 0 & 0 & 7 \\
 & 900 & 1618.9 & 14.4 & 0.2 & 1.5 & 5.0 & 1.2 & 60 & 100 & 0 & 33 & 0 & 0 \\
 & 1600 & 3598.0 & 114.4 & 0.7 & 0.8 & 17.6 & 9.3 & 50 & 100 & 0 & 100 & 0 & 7 \\\midrule
 
f & 900 & 83.3 & 4.2 & 5.5 & 2.0 & 3.7 & 173.3 & 0 & 82 & 0 & 0 & 0 & 100 \\
 & 1600 & 1319.5 & 17.5 & 10.5 & 9.6 & 24.2 & 737.8 & 0 & 77 & 0 & 0 & 0 & 100 \\
 & 3200 & 3580.7 & 200.7 & 20.3 & 2.2 & 164.2 & 3006.0 & 0 & 91 & 0 & 100 & 0 & 100 \\
 & 4000 & 3566.4 & 478.0 & 21.7 & 2.4 & 260.8 & 3569.0 & 0 & 84 & 0 & 100 & 0 & 100 \\
 & 5000 & 3551.2 & 1401.8 & 19.8 & 3.0 & 464.2 & 3580.6 & 0 & 84 & 0 & 100 & 23 & 100 \\\midrule
 
n & 100 & 0.2 & 0.1 & 0.0 & 0.1 & 0.2 & 0.1 & 55 & 100 & 30 & 0 & 0 & 29 \\
 & 200 & 3.7 & 0.5 & 0.0 & 0.4 & 1.1 & 0.4 & 43 & 100 & 16 & 0 & 0 & 26 \\
 & 400 & 28.6 & 1.9 & 0.1 & 1.2 & 2.9 & 1.9 & 58 & 100 & 8 & 0 & 0 & 18 \\
 & 900 & 1055.9 & 11.6 & 0.4 & 5.5 & 14.0 & 16.5 & 45 & 99 & 6 & 16 & 0 & 14 \\
 & 1600 & 2674.4 & 93.7 & 1.1 & 4.4 & 44.2 & 69.4 & 43 & 100 & 4 & 64 & 1 & 9 \\
 & 3200 & 3352.1 & 255.4 & 3.9 & 5.1 & 188.3 & 612.5 & 53 & 99 & 2 & 93 & 0 & 29 \\
 & 4000 & 3325.2 & 602.4 & 4.2 & 4.7 & 173.3 & 530.0 & 47 & 99 & 1 & 93 & 5 & 19 \\
 & 5000 & 3258.9 & 1348.0 & 5.2 & 10.3 & 280.0 & 716.0 & 41 & 99 & 2 & 92 & 25 & 27 \\
\end{tabular}
}
\caption{ISA: Statistics per class and number of jobs. Instances where LR terminated early without finding optimal are treated as timeout.}
\label{tab:isa-rmp_cg_times}
\end{table}

\section{Conclusion and Future Work}\label{sec:Conclusion}
It is worth emphasizing that Template pricing, while primarily discussed for set partitioning problems, can also be applied to other problem classes. The approach proves particularly beneficial when faced with highly degenerate cases. Our similarity function $d(x,y^i)$, was chosen for its simplicity, other choices could improve the results even further, especially on instances that have high degree of inherent fractionality at the root CG node.

Our Template pricing implementation efficiently uses a Lagrangian heuristic; however, other approaches may be more suitable when the subproblems are complex. Template pricing shares similarities to local search, that is, they both try to find suitable solutions without deviating too far from some baseline.  Local search around active columns is known to be a successful strategy~\cite{barnhart1998branch,vanderbeck_implementing_2005,vanderbeck_generic_2006}, typically with the assumption that checking valid perturbations (also known as exchange vectors) are much cheaper than solving the subproblem. While this is true, we have shown there can be additional benefits when solving the RMP. Our ideas could be easily incorporated into local search heuristics designed to approximate Template pricing for complex subproblems, for even further impact.

Another aspect that we have not addressed is that of a branch-and-price implementation. Finding ways to leverage templates within branching decisions presents an intriguing direction for future research.  Empirically, Template pricing is more likely to find integer solutions, so we imagine it could help reduce the branch-and-bound search.
There are several other interesting research areas to explore.  For example: combining Template pricing with other stabilization methods (like Pessoa); hybrid solutions with Lagrangian Relaxation providing initial dual bounds; modifications to the RMP (similar to DCA~\cite{elhallaoui2005dynamic}, or trust region hybrids).  

Lastly, we have focused our research on problem classes with independent subproblems, where the number of subproblems is known.  There are several problem classes that share subproblems across multiple columns, for example, bin packing or graph coloring.  In these problems, it is interesting to consider how a template is defined. We imagine that it could involve clustering the columns based on similarity within the current primal solution, or other heuristic approaches. Finally, a potential connection between recent theoretical results on escaping degeneracy for the simplex method, for example~\cite{kukharenko2024number}, and our proposed Template pricing approach could be an interesting direction to examine.

In closing, we have proposed Template pricing as a coordinated pricing paradigm for column generation that markedly accelerates convergence while incidentally finding high-quality integer solutions. We present both exact and heuristic variants, and show that their effectiveness stems from improved column compatibility and stability, which in turn speeds up the solution of the RMP. We further analyze the roles of initialization, column management, and degeneracy in overall performance. Computational results on GAP benchmarks show consistent and sometimes orders-of-magnitude improvements over Dantzig and Pessoa pricing.  Moreover, Template pricing is competitive with Lagrangian relaxation techniques, which have traditionally been the preferred approach for highly degenerate problems.


\bibliographystyle{splncs04}
\bibliography{bibliography}

\section{Statements and Declarations}
\subsection{Funding}
The authors declare that no funds, grants, or other support were received during the preparation of this manuscript.

\subsection{Competing Interests}
The authors have no relevant financial or non-financial interests to disclose.

\subsection{Author Contributions}
All authors contributed to the study conception and design. Implementation, material preparation, data collection and analysis were performed by Luke Marshall. An initial draft of the manuscript was prepared by Prachi Shah and substantially revised by Luke Marshall (incorporating additional algorithms and results), with input from all authors. All authors read and approved the final manuscript.

\subsection{Data Availability}

Datasets generated and analyzed during the current study are available in the \url{https://github.com/mathgeekcoder/template-pricing} repository.

\appendix

\section{Pessoa Directional Dual Smoothing}\label{sec:pessoa}

\begin{figure}[htbp]
\centering

\begin{tikzpicture}[
    scale=3,
    every node/.style={font=\small},
    state/.style={circle, fill=black, inner sep=1.6pt}
]
\coordinate (pi_hat) at (3,0);

\path (0,0) coordinate (pi_t)
      (0.5,1) coordinate (pi_g)
      (pi_t) -- (pi_g) coordinate[pos=0.8] (rho);

\path (pi_hat) -- (pi_t) coordinate[pos=0.475] (pi_k);
\path (pi_hat) -- (rho) coordinate[pos=0.475] (pi_tilde);
\path (pi_hat) -- (pi_g) coordinate[pos=0.35] (g_hat);

\def\offsetdist{0.5mm}

\draw[thick, Latex-Latex] 
    ($(pi_t)!\offsetdist!90:(rho)$) -- ($(rho)!\offsetdist!-90:(pi_t)$) 
    node[midway,sloped,above,allow upside down=false] {$\beta_t \norm{\pi^g_t - \pi_t}$};

\draw[thick, Latex-Latex] 
    ($(pi_k)!\offsetdist!90:(pi_t)$) -- ($(pi_t)!\offsetdist!-90:(pi_k)$) 
    node[midway,sloped,below,allow upside down=false] {$\alpha^k_t \norm{\hat{\pi}_{t-1} - \pi_t}$};

\draw[ultra thin] (pi_hat) -- (pi_t) -- (pi_g) -- (g_hat);
\draw[dashed, shorten >= 2pt] (pi_hat) -- (pi_tilde);
\draw[dashed, shorten <= 5pt] (pi_tilde) -- (rho);

\node[state,label=below left:$\pi_t$] at (pi_t) {};
\node[state,label=below right:$\hat{\pi}_{t-1}$] at (pi_hat) {};
\node[state,label=above left:$\pi^{g}_t$] at (pi_g) {};
\node[state,label=above left:$\rho_t$] at (rho) {};

\node[state,label=above left:$\pi^k_t$] at (pi_k) {};
\node[label={[shift={(-1.5mm,2mm)}]below left:$\tilde{\pi}_t$}] at (pi_tilde) {};

\path (pi_t) -- (pi_k)
        node[style={circle, fill=gray, inner sep=1pt}, pos=0.75, label={[gray]above:$\pi^2_t$}] {}
        node[style={circle, fill=gray, inner sep=1pt}, pos=0.5, label={[gray]above:$\pi^3_t$}] {}
        node[style={circle, fill=gray, inner sep=1pt}, pos=0.25, label={[gray]above:$\cdots$}] {}
        ;

\draw plot[mark=star, mark options={very thick, red, scale=0.5}] 
        coordinates {(pi_tilde)};

\draw[-Latex, thick, shorten >=3pt] (pi_k) -- (pi_tilde) {};

\draw[-{Latex[scale=1.5]}, shorten >= -1.1mm] (pi_hat) -- (g_hat);
\node[label={[shift={(1mm,-1mm)}]above right:$\hat{g}_{t-1}$}] at (g_hat) {};
\draw[thick] (pi_hat) -- (g_hat) {};

\pic [draw=black, thick, angle radius=2.3cm] {angle = pi_g--pi_hat--pi_t};

\path (pi_hat) -- (rho) coordinate[pos=0.125] (beta_label);
\node[label=left:$\cos^{-1} \beta_t$, yshift=-0.5mm] at (beta_label) {};

\end{tikzpicture}
\caption{Geometric view of Pessoa's directional smoothing and limited $\alpha^k_t$ backtracking line search. Terminates at the first point that yields a column with good reduced cost. Backtracking starts at *, when the gradient is available.}
\label{fig:pessoa}
\end{figure}
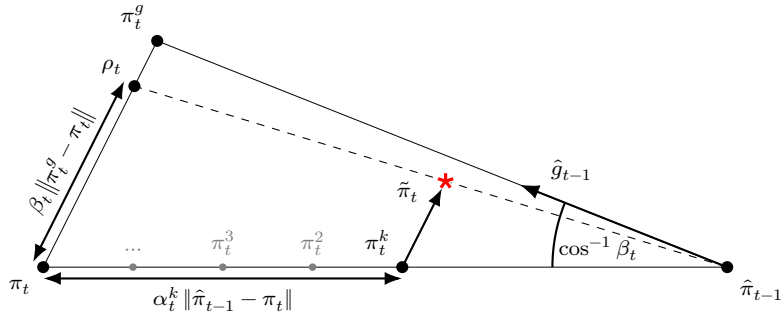

This section provides the full details of the directional dual smoothing algorithm defined by Pessoa, which is used in our computational study.  We include these details to clearly document the precise variant implemented, as Pessoa~\cite{pessoa_automation_2018} describe multiple stabilization schemes. See Figure \ref{fig:pessoa} for a geometric interpretation.

At each CG iteration $t$, we maintain the following information:
\begin{description}[leftmargin=!,labelwidth=\widthof{$\alpha_t$}]
    \item[$\pi_t$] dual vector obtained from solving the RMP
    \item[$\tilde{\pi}_t$] smoothed dual vector used for pricing
    \item[$\hat{\pi}_t$] ``best'' dual vector corresponding to the most recent iteration at which the RMP objective strictly improved
    \item[$\hat{g}_t$] subgradient associated with $\hat{\pi}_t$
    \item[$\alpha_t$] adaptive mixing parameter
    \item[$k$] backtracking index controlling smoothing intensity
\end{description}

For GAP, the pricing subproblems are decomposed by machine. Let $\{v^t_i\}_{i\in I}$ denote the set of columns produced by pricing at iteration $t$, one column for each machine $i\in I$.  Each subgradient element corresponding to job $j\in J$ is defined as: $$g_{tj} = \nabla\pi_{tj} = 1 - \sum_{i\in I} \mathbb{1}_{\{\mathbf{e}_j^\top v^t_{i} = 1\}}.$$
This is the same definition used by Lagrangian Relaxation in Section~\ref{sec:LR}. Intuitively, when:
\begin{description}[leftmargin=!,labelwidth=\widthof{$g_{tj}<0$}]
    \item[$g_{tj}=1$] job $j$ was not selected by any machine
    \item[$g_{tj}<0$] job $j$ was selected by multiple machines
    \item[$g_{t}=0$] the set of columns $\{v^t_i\}_{i\in I}$ partition the jobs across machines
\end{description}

Let $\hat{\pi}_{t-1}$ be the previous best dual vector and $\hat{g}_{t-1}$ be its associated subgradient.  Pessoa combines Wentges with an adaptive mixing and a directional step along the subgradient.  First, consider the family of convex combinations, corresponding to the backtracking index $k=1, \dots, 9$ and the previous adaptive mixing parameter $\alpha_{t-1}$:
\begin{align*}
\alpha^k_t &= \max \left\{0,\; 1 - k\left(1-\alpha_{t-1}\right)\right\} \\
\pi^k_t &= \alpha^k_t\hat{\pi}_{t-1} + \left(1-\alpha^k_t\right)\pi_t
\end{align*}
Next, if $\norm{\hat{g}_t}>0$, the directional step is calculated as:
\begin{align*}
\pi^g_t &= \hat{\pi}_{t-1} + \norm{\pi_t - \hat{\pi}_{t-1}}\frac{\hat{g}_{t-1}}{\norm{\hat{g}_{t-1}}}\\
\hspace{0.1\textwidth}\beta_t &= \frac{(\pi_t - \hat{\pi}_{t-1})\cdot(\pi^g_t - \hat{\pi}_{t-1})}{\norm{\pi_t - \hat{\pi}_{t-1}} \norm{\pi^g_t - \hat{\pi}_{t-1}}} \\
\rho_t &= \beta_t\pi^g_t + (1-\beta_t)\pi_t \\
\end{align*}
such that the smoothed dual vector is given by:
$$\tilde{\pi}_t =
\max \left\{\mathbf{0},\; \hat{\pi}_{t-1} + \norm{\pi^k_t - \hat{\pi}_{t-1}} \frac{\rho_t - \hat{\pi}_{t-1}}{\norm{\rho_t - \hat{\pi}_{t-1}}} \right\}
$$
Non-negativity is enforced in the definition above, that is, $\max \left\{ \mathbf{0},\; \dots\right\}$ is taken element wise.  

Pessoa performs a limited backtracking search to guarantee finding at least one column having negative reduced cost (if one exists). That is, minimizing $k$ such that there exists some $i\in I$ having 
\begin{align*}
    x^i = \underset{\mathbf{x}\in P_i}{\arg\min} &\sum_{j \in J} (c_{ij} - \tilde{\pi}_{tj}) x_j, \\
    \text{with } &\sum_{j \in J} \left(c_{ij} - \pi_{tj}\right) x^i_j \le \mu_{ti} - \varepsilon.
\end{align*}
However, if we require any backtracking, that is, $k>1$, or when $\norm{\hat{g}_t}=0$, the smoothed dual vector is simply given by:
$$\tilde{\pi}_t = \pi^k_t$$
Notice that as $k$ increases (more backtracking), $\alpha^k_t\rightarrow 0$ and $\tilde{\pi}_t\rightarrow\pi_t$. If no negative reduced columns can be found for all $k\in\{1, \dots, 9\}$, Pessoa reverts to Dantzig pricing:
$$\tilde{\pi}_t = \pi_t$$
The final step is to update the adaptive mixing parameter for the next iteration:
$$\alpha_t = \begin{cases}
\min\left\{0.9999,\; 0.9\alpha_{t-1} + 0.1\right\} & \quad\text{if }\, \tilde{g}_t \cdot (\pi_t - \hat{\pi}_{t-1}) > 0  \\
\max\left\{0,\; \alpha_{t-1} - 0.1\right\} & \quad \text{o/w}
\end{cases}$$

\section{Solving Template pricing exactly via hierarchical optimization}\label{sec:hie}
We can solve the Template pricing problem exactly using a hierarchical multiple objective integer program.  
Formally, here is the sequence of optimization problems we solve:
\begin{enumerate}
    \item Verify $S^i_{(\pi,\mu)} \neq \emptyset$, that is, solve the Dantzig pricing problem (\ref{eq:dantzig-pricing}). If non-negative reduced costs, then we have $S^i_{(\pi,\mu)} = \emptyset$ and we stop.
    
    \item Else, optimize template objective over the \textit{good} columns, that is, solve~\eqref{eq:objinner}:
    \begin{eqnarray*}
    \textup{OPT}:= &\max& d\left(x,y^i\right) \\
    &\textup{s.t.}& x \in S^i_{(\pi,\mu)} \equiv (\ref{eq:goodredcost})
    \end{eqnarray*}
    
    \item Finally, break ties among the optimal solutions of (\ref{eq:objinner}), that is, solve~\eqref{eq:objouter}:
    \begin{eqnarray}\label{eq:complex}
    \begin{array}{rcl}
    &\min& \displaystyle \sum_{j \in J} \left(c_{ij} - \pi_j\right)x_j\\
    &\textup{s.t.}& d(x,y^i) \geq \textup{OPT} \\
    && x\in P_i
    \end{array}
    \end{eqnarray}
\end{enumerate}
Due to numerical issues, \eqref{eq:complex} can fail at times to find a column with negative reduced costs, even when one exists. In these cases, we reduce the value of OPT by $1$ and re-solve (\ref{eq:complex}) until a suitable column is found.

\section{Column management policy}
\label{app:column-management}

This section describes how the column management policy mentioned in Section~\ref{sec:column-management} was constructed.

Recall that our implementation uses age-based retention. For each column, we record its age, equal to the most recent iteration in which it was in the optimal RMP basis; or the iteration when it was first added to the RMP. At iteration $t$, columns are retained when their recorded age is in the interval $\{t-\tau, \ldots, t\}$, where $\tau$ denotes the age threshold.

\begin{figure}[htbp]
\centering
\colorlet{mylow}{green}
\colorlet{mymed}{Goldenrod}
\colorlet{myhigh}{msred}

\resizebox{0.8\columnwidth}{!}{%
\begin{tikzpicture}[scale=0.3]
\tikzset{mysquare/.style={minimum size=3.0mm, draw, inner sep=0}}
  \pgfplotstableread[col sep=tab,header=false]{c10400-heat-v2.tsv}{\datatable}

  \pgfplotstablegetcolsof{\datatable}
  \pgfmathtruncatemacro{\colcount}{\pgfplotsretval-1}
  \node[] at (\colcount*0.5, 1.5) {\tiny Age-Threshold (\# iterations)};
  \node[] at (28cm, 1.5) {\tiny Sparklines};
  \node[rotate=90] at (-3, -5.5) {\tiny Activation Threshold};
  
  \foreach \col in {1, 5, 10, ..., \colcount} {%
    \pgfplotstablegetelem{0}{[index]\col}\of\datatable
    \edef\colname{\pgfplotsretval}
    \node[minimum size=0mm] at (\col,0.0) {\tiny \colname};
  }
 
  \pgfplotstablegetrowsof{\datatable}
  \pgfmathtruncatemacro{\rows}{\pgfplotsretval-1}
  \foreach \r in {1,...,10} {
    \foreach \c in {1,...,\colcount} { 
      \pgfplotstablegetelem{\r}{[index]\c}\of\datatable
      \edef\cell{\pgfplotsretval}
      \ifx\cell\empty\else
        \pgfmathsetmacro{\val}{\cell}
        \pgfmathtruncatemacro{\pct}{(100*\val)}

    \begingroup
        \pgfmathsetmacro{\threshold}{20}%
        \ifdim \pct bp < \threshold bp%
            \colorlet{c1}{mylow}
            \colorlet{c2}{mymed}
            \pgfmathsetmacro{\mix}{(100 * (\pct / \threshold)}%
        \else%
            \pgfmathsetmacro{\mix}{min(100, 100 * ((\pct - \threshold) / (50 - \threshold)))}%
            \colorlet{c1}{mymed}
            \colorlet{c2}{myhigh}
        \fi
    
        \colorlet{testcolor}{c2!\mix!c1}
        \node[mysquare, draw=none,line width=0.00mm,fill=testcolor, rectangle] at (\c,-\r) {};
    \endgroup
      \fi
    }
    
\begingroup
      \pgfmathsetmacro{\rof}{\r+0.48}
\begin{axis}[
  axis lines=none,
  ticks=none,
  enlargelimits=0.05,
  width=6.5cm,
  height=25.5mm,
  shift={(25.5cm, -\rof cm)},
]
\addplot [very thick] table[x index=0, y index=\r] {c10400-heat-v2-transpose.tsv};
\end{axis};

\endgroup
    \pgfplotstablegetelem{\r}{[index]0}\of\datatable
    \edef\rl{\pgfplotsretval}
    \node[text width=3mm, align=right] at (-0.5,-\r) {\tiny \rl};
  }
\end{tikzpicture}
}
\caption{Column management heatmap for Dantzig pricing on instance c10400. Solve time is shown as a function of age threshold and activation frequency. Green and yellow correspond to shorter solve times, while red indicates 1-hour timeout.  Each cell is a locally smoothed average of 200 replications including nearby parameter settings. The sparklines on the right show the relative solve-time profile for each activation frequency setting.} \label{fig:heatmap}
\end{figure}
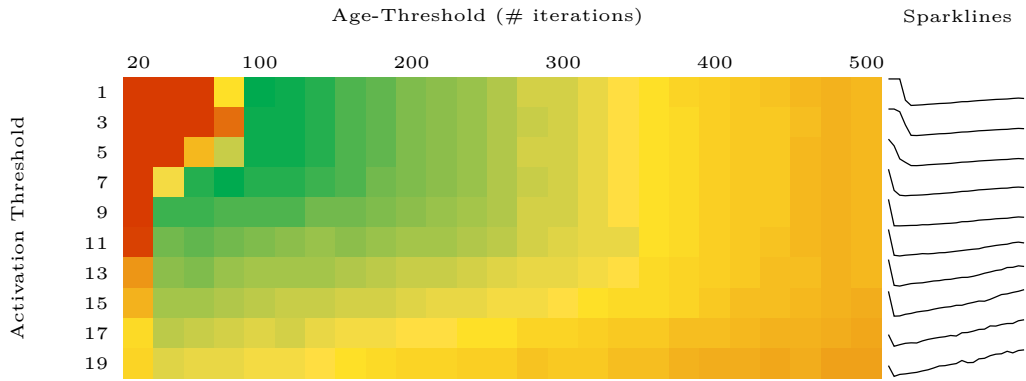

In addition to $\tau$, we also considered an activation-frequency parameter controlling when column removal is triggered. This parameter was implemented as a multiplier on the number of rows in the RMP, following guidance from the literature~\cite{desrosiers_branch-and-price_2025}, so that activation scales appropriately with problem size. We first performed a two-parameter sweep over age threshold and activation frequency on representative instances. Figure~\ref{fig:heatmap} shows an example. These preliminary experiments indicate that, in our setting, activating column removal at every iteration was sufficient to minimize solve times. Since the RMP is solved via primal simplex, this choice has little downside in practice because the basis itself remains unchanged. Accordingly, we fixed activation at every iteration in the following experiments and those reported in the main text.

The choice for age threshold $\tau$ is a tradeoff. Small values of $\tau$ can occasionally yield very fast solves when the retained columns happen, by chance, to be sufficient; however, this behavior is unreliable and typically leads to volatile performance. Larger values of $\tau$ are more likely to retain columns that provide a good approximation of the dual polyhedron; improving robustness, at the cost of larger RMPs and increased solve times.

\begin{figure}[htbp]
\centering
\begin{tikzpicture}
\begin{axis}[
  name=left fig,
  axis x line*=bottom,
  xmin=0, xmax=250,
  ymin=0, ymax=1004,
  ylabel={Time (s)},
  ylabel near ticks,
  xlabel={Age threshold (\# its)},
  tick label style={font=\small}, 
  width=0.95\columnwidth,
  height=4.5cm,
  legend cell align={left},
  legend style={/tikz/every even column/.append style={column sep=0.5cm}, font = \small, at={(0.5,1.3)},legend cell align={left}, anchor=north,legend columns=-1, draw=none, fill=none}
]

\addplot [draw=msblue, only marks, mark=*, mark options={fill=msblue,solid,fill opacity=0.1,draw opacity=0.5, scale=0.5}] 
    table[x=age_limit, y=time] {a05200-dantzig-age-sweep.tsv};

\addplot [draw=Black, only marks, mark=*, mark options={fill=Black,solid,fill opacity=0.9,draw opacity=0.5, scale=0.5}] 
    table[x=age_limit, y=time] {a05200-dantzig-age-sweep-smoothed.tsv};

\draw [dotted, semithick] (axis cs: 132, 0) -- (axis cs: 132, 1004);
\end{axis}
\end{tikzpicture} 
\caption{Age-threshold sweep for Dantzig pricing on instance a05200. Solve time is shown as a function of the age threshold using five replications and a 1000-second time limit. Black solid dots show the centered rolling geometric mean (window size 5), and the dashed vertical line marks the selected threshold. The figure illustrates the high variability at small thresholds and the robustness-based criteria used to select the threshold.}
\label{fig:a05200-retention-sweep}
\end{figure}
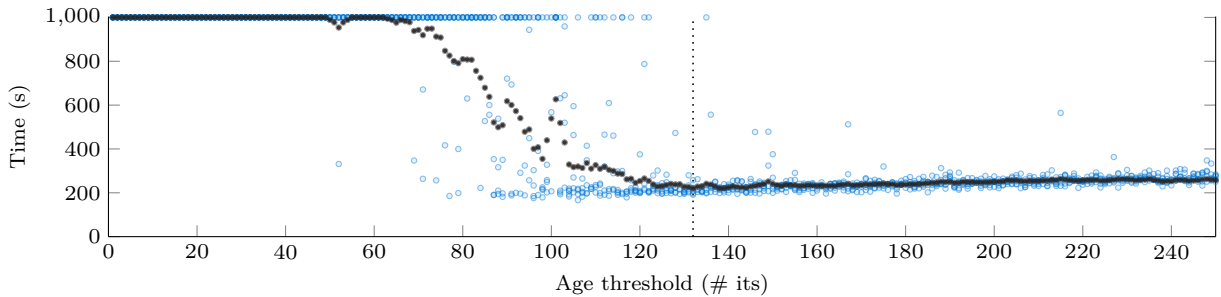

To calibrate $\tau$, we performed parameter sweeps on a representative set of instances. For each instance and candidate $\tau$, we ran five replications and measured time to convergence. To reduce noise, we computed the geometric mean across replications and then applied a centered rolling weighted average (across five data points). We selected the smallest threshold whose smoothed value was effectively tied with the minimum, where ties were defined as being within 1\% or 1 second of the best observed value. Figure~\ref{fig:a05200-retention-sweep} illustrates the selection procedure on a single instance.

We repeated this procedure for Dantzig, Pessoa, and LT pricing. Our assumption is that the best age-threshold policy is driven by instance degeneracy, for which we use the job-to-machine ratio as a proxy. 
To define our policy, rather than fitting a polynomial to the selected thresholds, we construct a quadratic \emph{covering} policy. Specifically, for each method, we collect the selected threshold from every instance, and then solve a quadratic program to find polynomial coefficients such that the resulting curve covers these data points while remaining as small as possible. We use a cover rather than a fit because a fitted curve could assign $\tau$ values below the empirically selected ``best'' thresholds, thereby increasing the risk of choosing values that are too aggressive and hence too volatile. Requiring the policy to lie above the selected thresholds yields a more conservative and robust rule. 

Figure~\ref{fig:retention-sweep} shows the chosen policy and the selected thresholds (per instance) against the job-to-machine ratio. The chosen thresholds increase with degeneracy, and are well approximated by a quadratic function. These quadratic covering curves define the age-threshold policy used in the experiments in Section~\ref{sec:results}.


\section{Full Yagiura results}
\begin{table}[htbp]
\centering
\newcolumntype{d}{R{14mm}{21600}{}{2.5mm}}
\newcolumntype{g}{R{14mm}{4000}{}{2.5mm}}
\newcolumntype{t}{R{14mm}{24410}{Dandelion}{2.5mm}}
\resizebox{0.81\columnwidth}{!}{%
\addtolength{\tabcolsep}{-0.4em}
\begin{tabular}{rrrddddggggtt}
& && \multicolumn{4}{c}{RMP (s)} & \multicolumn{4}{c}{Pricing (s)} & \multicolumn{2}{c}{$\times$ Time vs LT} \\ \cmidrule(l){4-7}\cmidrule(l){8-11}\cmidrule(l){12-13}
class & $|J|$ & $|I|$ & D & P & LT & MT & D & P & LT & MT & D & P \EndTableHeader\\ \toprule

a & 100  &  5 & 2.0 & 0.2 & 0.0 & 0.0 & 0.2 & 0.1 & 0.0 & 0.0 & 452 & 59 \\
 &       & 10 & 0.5 & 0.1 & 0.0 & 0.0 & 0.1 & 0.0 & 0.0 & 1.4 &  18 &  6 \\
 &       & 20 & 0.2 & 0.1 & 0.0 & 0.0 & 0.0 & 0.0 & 0.0 & 0.0 &  44 & 20\\
\cmidrule{3-13}
 & 200 &  5 & 117.9 & 5.0 & 0.0 & 0.0 & 0.9 & 0.2 & 0.0 & 0.0 & 23224 & 1022 \\
 &     & 10 &  13.7 & 1.2 & 0.0 & 0.0 & 0.2 & 0.1 & 0.0 & 0.1 &  2380 &  211 \\
 &     & 20 &   3.3 & 1.2 & 0.0 & 0.0 & 0.1 & 0.1 & 0.0 & 1.6 &    64 &   24 \\
\midrule
b & 100 &  5 & 1.5 & 0.3 & 0.1 & 0.2 & 0.2 & 0.1 & 0.1 & 37.7 & 11 & 3 \\
 &      & 10 & 0.4 & 0.1 & 0.0 & 0.1 & 0.1 & 0.0 & 0.0 &  3.7 &  9 & 3 \\
 &      & 20 & 0.2 & 0.1 & 0.0 & 0.0 & 0.0 & 0.0 & 0.0 &  1.9 &  7 & 4 \\
\cmidrule{3-13}
 & 200 &  5 & 96.3 & 6.6 & 0.8 & 1.2 & 1.0 & 0.4 & 0.5 & 106.5 & 75 & 5 \\
 &     & 10 & 12.4 & 1.4 & 0.2 & 0.3 & 0.2 & 0.1 & 0.1 &  24.1 & 43 & 5 \\
 &     & 20 &  2.6 & 0.6 & 0.1 & 0.1 & 0.1 & 0.0 & 0.0 &   7.5 & 26 & 6 \\
\midrule

c & 100 &  5 & 1.6 & 0.3 & 0.1 & 0.2 & 0.1 & 0.1 & 0.1 & 22.4 & 11 & 3 \\
 &      & 10 & 0.4 & 0.1 & 0.0 & 0.1 & 0.0 & 0.0 & 0.0 &  5.4 &  7 & 3 \\
 &      & 20 & 0.1 & 0.1 & 0.0 & 0.0 & 0.0 & 0.0 & 0.0 &  2.0 &  5 & 3 \\
\cmidrule{3-13}
 & 200 &  5 & 96.6 & 6.2 & 0.5 & 0.8 & 1.1 & 0.4 & 0.4 & 88.4 & 100 & 7 \\
 &     & 10 & 11.5 & 1.4 & 0.2 & 0.3 & 0.2 & 0.1 & 0.1 & 25.4 &  37 & 5 \\
 &     & 20 &  2.5 & 0.6 & 0.1 & 0.1 & 0.1 & 0.0 & 0.0 &  8.4 &  21 & 5 \\
 \cmidrule{3-13}
 & 400 & 10 & 932.3 & 44.5 & 1.6 & 3.8 & 1.4 & 0.5 & 0.9 & 196.3 & 381 & 18 \\
 &     & 20 & 104.4 &  9.5 & 0.5 & 1.4 & 0.2 & 0.2 & 0.2 &  92.6 & 144 & 13 \\
 &     & 40 &  22.2 &  3.7 & 0.2 & 0.4 & 0.1 & 0.1 & 0.1 &  11.6 &  62 & 11 \\
\cmidrule{3-13}
 & 900 & 15 & 21598.0 & 2998.3 & 11.1 & 19.0 & 1.3 & 2.9 & 2.8 & 263.5 & 1538 & 214 \\
 &     & 30 &  7507.9 &  336.1 &  3.9 &  6.3 & 1.0 & 0.5 & 0.7 & 136.3 & 1587 &  71 \\
 &     & 60 &  1054.6 &   84.1 &  2.1 &  2.6 & 0.2 & 0.2 & 0.3 &  31.4 &  402 &  32 \\
\cmidrule{3-13}
 & 1600 & 20 & 21600.0 & 21598.7 & 90.8 & 81.7 & 0.8 & 1.4 & 8.7 & 677.0 &  221 & 221 \\
 &      & 40 & 21600.0 &  7752.2 & 25.5 & 20.5 & 0.4 & 1.4 & 2.0 & 244.1 &  786 & 282 \\
 &      & 80 & 21600.0 &  1515.9 & 11.2 &  9.1 & 0.3 & 0.5 & 0.9 &  61.4 & 1695 & 118 \\
 \midrule
d & 100 &  5 & 1.6 & 0.4 & 0.2 & 0.3 & 0.2 & 0.2 & 0.4 & 81.5 & 3 & 1 \\
 &      & 10 & 0.4 & 0.2 & 0.1 & 0.1 & 0.1 & 0.1 & 0.1 & 31.4 & 2 & 1 \\
 &      & 20 & 0.2 & 0.1 & 0.1 & 0.1 & 0.0 & 0.0 & 0.1 & 17.6 & 2 & 1 \\ 
\cmidrule{3-13}
 & 200 &  5 & 85.1 & 8.4 & 1.9 & 1.5 & 2.8 & 1.3 & 3.1 & 251.1 & 18 & 2 \\
 &     & 10 & 13.1 & 2.6 & 0.9 & 0.7 & 0.4 & 0.3 & 1.3 & 143.7 &  6 & 1 \\
 &     & 20 &  3.5 & 1.2 & 0.5 & 0.5 & 0.1 & 0.1 & 0.4 &  63.7 &  4 & 1 \\
 \cmidrule{3-13}
 & 400 & 10 & 1060.6 & 85.7 & 6.6 & 4.2 & 7.7 & 3.4 & 8.5 & 317.1 & 70 & 6 \\
 &     & 20 &  160.6 & 25.5 & 3.0 & 2.4 & 1.2 & 0.8 & 3.2 & 293.3 & 26 & 4 \\
 &     & 40 &   44.6 & 11.3 & 1.7 & 1.9 & 0.3 & 0.3 & 1.0 &  95.3 & 16 & 4 \\
\cmidrule{3-13}
 & 900 & 15 & 21594.0 & 5043.5 & 33.9 & 33.5 &  5.7 & 25.9 & 36.4 & 1033.3 & 307 & 72 \\
 &     & 30 & 14729.0 & 1271.3 & 15.6 & 15.4 & 10.1 &  6.5 & 16.4 &  517.3 & 454 & 39 \\
 &     & 60 &  3370.0 &  500.4 &  7.4 &  6.4 &  2.1 &  1.8 &  4.9 &  137.9 & 265 & 39 \\
 \cmidrule{3-13}
 & 1600 & 20 & 21600.0 & 21594.5 & 233.7 & 191.4 & 2.6 &  7.2 & 106.3 & 1569.1 &  64 &  64 \\
 &      & 40 & 21596.0 & 21594.8 &  51.0 &  37.8 & 3.1 &  5.9 &  32.7 &  405.1 & 256 & 256 \\
 &      & 80 & 21599.0 & 11758.9 &  26.7 &  19.3 & 0.8 &  5.9 &  20.8 &  268.3 & 441 & 238 \\
\midrule
e & 100 &  5 & 1.4 & 0.3 & 0.1 & 0.2 & 0.1 & 0.1 & 0.1 & 43.4 & 9 & 3 \\
 &      & 10 & 0.4 & 0.2 & 0.1 & 0.1 & 0.1 & 0.0 & 0.0 & 15.8 & 4 & 2 \\
 &      & 20 & 0.2 & 0.2 & 0.1 & 0.1 & 0.0 & 0.0 & 0.0 & 12.3 & 3 & 2 \\
\cmidrule{3-13}
 & 200 &  5 & 86.4 & 7.4 & 0.9 & 1.3 & 0.9 & 0.4 & 0.6 & 97.2 & 62 & 6 \\
 &     & 10 & 11.9 & 2.6 & 0.3 & 0.6 & 0.2 & 0.1 & 0.3 & 34.6 & 20 & 4 \\
 &     & 20 &  3.2 & 1.4 & 0.2 & 0.5 & 0.1 & 0.1 & 0.1 & 46.3 & 10 & 4 \\
 \cmidrule{3-13}
 & 400 & 10 & 878.4 & 70.9 & 3.3 & 4.5 & 1.8 & 0.9 & 1.6 & 155.6 & 181 & 15 \\
 &     & 20 & 129.6 & 22.8 & 1.4 & 2.6 & 0.4 & 0.3 & 0.8 & 208.9 &  57 & 10 \\
 &     & 40 &  30.1 & 11.4 & 0.8 & 1.7 & 0.1 & 0.1 & 0.2 & 159.3 &  24 &  9 \\
\cmidrule{3-13}
 & 900 & 15 & 21597.0 & 4400.4 & 34.0 & 74.1 & 2.0 & 6.5 & 9.3 & 2244.2 & 543 & 111 \\
 &     & 30 & 10069.0 &  986.8 & 10.0 & 16.6 & 2.6 & 1.9 & 3.7 & 1146.5 & 713 &  70 \\
 &     & 60 &  1414.9 &  323.6 &  7.0 & 12.6 & 0.6 & 0.6 & 1.5 & 1193.1 & 150 &  34 \\
\cmidrule{3-13}
 & 1600 & 20 & 21599.0 & 21599.3 & 167.3 & 160.6 & 0.7 & 2.0 & 29.2 & 1950.1 & 111 & 111 \\
 &      & 40 & 21599.0 & 21017.8 &  63.7 &  70.5 & 0.6 & 8.4 & 12.3 & 3067.6 & 279 & 272 \\
 &      & 80 & 21596.0 &  5896.3 &  59.4 &  37.3 & 0.6 & 3.0 &  8.8 & 3596.6 & 312 &  85 \\  
\end{tabular}
}
\caption{Average of RMP time (columns 4–7) and Pricing time (columns 8–11) across classes, jobs, and machines. Additionally, total time (RMP + Pricing) comparison against LT (columns 12-13), clearly showing instances where Template pricing is over 1000$\times$ faster than Dantzig, and 100$\times$ faster than Pessoa. Five replications each, six hour time limit.}
\label{tab:yagiura-full}
\end{table}

\end{document}